# BOOSTING WITH EARLY STOPPING: CONVERGENCE AND CONSISTENCY

By Tong Zhang and Bin Yu[1]

*IBM T. J. Watson Research Center and University of California at Berkeley*

Boosting is one of the most significant advances in machine learning for classification and regression. In its original and computationally flexible version, boosting seeks to minimize empirically a loss function in a greedy fashion. The resulting estimator takes an additive function form and is built iteratively by applying a base estimator (or learner) to updated samples depending on the previous iterations. An unusual regularization technique, early stopping, is employed based on CV or a test set.

This paper studies numerical convergence, consistency and statistical rates of convergence of boosting with early stopping, when it is carried out over the linear span of a family of basis functions. For general loss functions, we prove the convergence of boosting's greedy optimization to the infimum of the loss function over the linear span. Using the numerical convergence result, we find early-stopping strategies under which boosting is shown to be consistent based on i.i.d. samples, and we obtain bounds on the rates of convergence for boosting estimators. Simulation studies are also presented to illustrate the relevance of our theoretical results for providing insights to practical aspects of boosting.

As a side product, these results also reveal the importance of restricting the greedy search step-sizes, as known in practice through the work of Friedman and others. Moreover, our results lead to a rigorous proof that for a linearly separable problem, AdaBoost with $\varepsilon \to 0$ step-size becomes an $L^1$-margin maximizer when left to run to convergence.

**1. Introduction.** In this paper we consider boosting algorithms for classification and regression. These algorithms represent one of the major ad-

---

Received January 2003; revised September 2004.

[1]Supported in part by NSF Grant FD01-12731 and ARO Grant DAAD19-01-1-0643.

vances in machine learning. In their original version, the computational aspect is explicitly specified as part of the estimator/algorithm. That is, the empirical minimization of an appropriate loss function is carried out in a greedy fashion, which means that at each step a basis function that leads to the largest reduction of empirical risk is added into the estimator. This specification distinguishes boosting from other statistical procedures which are defined by an empirical minimization of a loss function without the numerical optimization details.

Boosting algorithms construct composite estimators using often simple base estimators through the greedy fitting procedure. An unusual regularization technique, early stopping, is employed based on CV or a test set. This family of algorithms has been known as the stagewise fitting of *additive models* in the statistics literature [18, 17]. For the squared loss function, they were often referred to in the signal processing community as *matching pursuit* [29]. More recently, it was noticed that the AdaBoost method proposed in the machine learning community [13] can also be regarded as stagewise fitting of additive models under an exponential loss function [7, 8, 15, 31, 34]. In this paper we use the term *boosting* to indicate a greedy stagewise procedure to minimize a certain loss function empirically. The abstract formulation will be presented in Section 2.

Boosting procedures have drawn much attention in the machine learning community as well as in the statistics community, due to their superior empirical performance for classification problems. In fact, boosted decision trees are generally regarded as the best off-the-shelf classification algorithms we have today. In spite of the significant practical interest in boosting, a number of theoretical issues have not been fully addressed in the literature. In this paper we hope to fill some gaps by addressing three basic issues regarding boosting: its numerical convergence when the greedy iteration increases, in Section 4.1; its consistency (after early stopping) when the training sample size gets large, in Sections 3.3 and 5.2; and bounds on the rate of convergence for boosting estimators, in Sections 3.3 and 5.3.

It is now well known that boosting forever can overfit the data (e.g., see [16, 19]). Therefore, in order to achieve consistency, it is necessary to stop the boosting procedure early (but not too early) to avoid overfitting. In the early stopping framework, the consistency of boosting procedures has been considered by Jiang for exponential loss [19] boosting (but the consistency is in terms of the classification loss) and Bühlmann under squared loss [10] for tree-type base classifiers. Jiang's approach also requires some smoothness conditions on the underlying distribution, and it is nonconstructive (hence does not lead to an implementable early-stopping strategy). In Sections 3.3 and 5.2 we present an early-stopping strategy for general loss functions that guarantees consistency.

A different method of achieving consistency (and obtaining rate of convergence results) is through restricting the weights of the composite estimator using the 1-norm of its coefficients (with respect to the basis functions). For example, this point of view is taken up in [5, 28, 30]. In this framework, early stopping is not necessary since the degree of overfitting or regularization is controlled by the 1-norm of the weights of the composite estimator. Although this approach simplifies the theoretical analysis, it also introduces an additional control quantity which needs to be adjusted based on the data. Therefore, in order to select an optimal regularization parameter, one has to solve many different optimization problems, each with a regularization parameter. Moreover, if there are an infinite (or extremely large) number of basis functions, then it is not possible to solve the associated 1-norm regularization problem. Note that in this case greedy boosting (with approximate optimization) can still be applied.

A question related to consistency and rate of convergence is the convergence of the boosting procedure as an optimization method. This is clearly one of the most fundamental theoretical issues for boosting algorithms. Previous studies have focused on special loss functions. Specifically, Mallat and Zhang proved the convergence of matching pursuit in [29], which was then used in [10] to study consistency; in [9] Breiman obtained an infinite-sample convergence result of boosting with the exponential loss function for ±1-trees (under some smoothness assumptions on the underlying distribution), and the result was used by Jiang to study the consistency of AdaBoost. In [12] a Bregman divergence-based analysis was given. A convergence result was also obtained in [31] for a gradient descent version of boosting.

None of these studies provides any information on the numerical speed of convergence for the optimization. The question of numerical speed of convergence has been studied when one works with the 1-norm regularized version of boosting where we assume that the optimization is performed in the convex hull of the basis functions. Specifically, for function estimation under least-squares loss, the convergence of the greedy algorithm in the convex hull was studied in [1, 20, 25]. For general loss functions, the convergence of greedy algorithms (again, the optimization is restricted to the convex hull) was recently studied in [37]. In this paper we apply the same underlying idea to the standard boosting procedure where we do not limit the optimization to the convex hull of the basis functions. The resulting bound provides information on the speed of convergence for the optimization. An interesting observation of our analysis is the important role of small step-size in the convergence of boosting procedures. This provides some theoretical justification for Friedman's empirical observation [14] that using small step-sizes almost always helps in boosting procedures.

Moreover, the combination of numerical convergence results with modern empirical process bounds (based on Rademacher complexity) provides a way to derive bounds on the convergence rates of early-stopping boosting procedures. These results can be found in Sections 3.3 and 5.3. Section 6 contains a simulation study to show the usefulness of the insights from our theoretical analyses in practical implementations of boosting. The proofs of the two main results in the numerical convergence section (Section 4.1) are deferred to Section A.2. Section A.3 discusses relaxations of the restricted step-size condition used for earlier results, and Section A.4 uses numerical convergence results to give a rigorous proof of the fact that for separable problems, AdaBoost with small step-size becomes an $L_1$ margin maximizer at its limit (see [18]).

**2. Abstract boosting procedure.** We now describe the basics to define the boosting procedure that we will analyze in this paper. A similar setup can be found in [31]. The main difference is that the authors in [31] use a gradient descent rule in their boosting procedure while here we use approximate minimization.

Let $S$ be a set of real-valued functions and define

$$\text{span}(S) = \left\{ \sum_{j=1}^{m} w^j f^j : f^j \in S, w^j \in R, m \in Z^+ \right\},$$

which forms a linear function space. For all $f \in \text{span}(S)$, we can define the 1-norm with respect to the basis $S$ as

$$(1) \qquad \|f\|_1 = \inf \left\{ \|w\|_1; f = \sum_{j=1}^{m} w^j f^j : f^j \in S, m \in Z^+ \right\}.$$

We want to find a function $\bar{f} \in \text{span}(S)$ that approximately solves the optimization problem

$$(2) \qquad \inf_{f \in \text{span}(S)} A(f),$$

where $A$ is a convex function of $f$ defined on $\text{span}(S)$. Note that the optimal value may not be achieved by any $f \in \text{span}(S)$, and for certain formulations (such as AdaBoost) it is possible that the optimal value is not finite. Both cases are still covered by our results, however.

The abstract form of the greedy boosting procedure (with restricted step-size) considered in this paper is given by the following algorithm:

ALGORITHM 2.1 (Greedy boosting).

Pick $f_0 \in \text{span}(S)$
**for** $k = 0, 1, 2, \ldots$
    Select a closed subset $\Lambda_k \subset R$ such that $0 \in \Lambda_k$ and $\Lambda_k = -\Lambda_k$
    Find $\bar{\alpha}_k \in \Lambda_k$ and $\bar{g}_k \in S$ to approximately minimize the function:

$(*)$                                            $(\alpha_k, g_k) \to A(f_k + \alpha_k g_k)$

Let $f_{k+1} = f_k + \bar{\alpha}_k \bar{g}_k$
**end**

REMARK 2.1.   The approximate minimization of $(*)$ in Algorithm 2.1 should be interpreted as finding $\bar{\alpha}_k \in \Lambda_k$ and $\bar{g}_k \in S$ such that

$$(3) \qquad A(f_k + \bar{\alpha}_k \bar{g}_k) \leq \inf_{\alpha_k \in \Lambda_k, g_k \in S} A(f_k + \alpha_k g_k) + \varepsilon_k,$$

where $\varepsilon_k \geq 0$ is a sequence of nonnegative numbers that converges to zero.

REMARK 2.2.   The requirement that $0 \in \Lambda_k$ is not crucial in our analysis. It is used as a convenient assumption in the proof of Lemma 4.1 to simplify the conditions. Our convergence analysis allows the choice of $\Lambda_k$ to depend on the previous steps of the algorithm. However, the most interesting $\Lambda_k$ for the purpose of this paper will be independent of previous steps of the algorithm:

  (a)  $\Lambda_k = R$,
  (b)  $\sup \Lambda_k = \tilde{h}_k$ where $\tilde{h}_k \geq 0$ and $\tilde{h}_k \to 0$.

As we will see later, the restriction of $\alpha_k$ to the subset $\Lambda_k \subset R$ is useful in the convergence analysis.

As we shall see later, the step-size $\bar{\alpha}_k$ plays an important role in our analysis. A particular interesting case is to restrict the step-size explicitly. That is, we assume that the starting point $f_0$, as well as quantities $\varepsilon_k$ and $\Lambda_k$ in (3), are sample-independent, and $h_k = \sup \Lambda_k$ satisfies the conditions

$$(4) \qquad \sum_{j=0}^{\infty} h_j = \infty, \qquad \sum_{j=0}^{\infty} h_j^2 < \infty.$$

The reason for this condition will become clear in the numerical convergence analysis of Section 4.1.

**3. Assumptions and main statistical results.**   The purpose of this section is to state assumptions needed for the analyses to follow, as well as the main statistical results. There are two main aspects of our analysis. The first is the numerical convergence of the boosting algorithm as the number of iterations increases, and the second is the statistical convergence of the resulting boosting estimator, so as to avoid overfitting. We list respective assumptions separately. The statistical consistency result can be obtained by combining these two aspects.

3.1. *Assumptions for the numerical convergence analysis.* For all $f \in$ span$(S)$ and $g \in S$, we define a real-valued function $A_{f,g}(\cdot)$ as

$$A_{f,g}(h) = A(f + hg).$$

DEFINITION 3.1.   Let $A(f)$ be a function of $f$ defined on span$(S)$. Denote by span$(S)'$ the dual space of span$(S)$ [i.e., the space of real-valued linear functionals on span$(S)$]. We say that $A$ is differentiable with gradient $\nabla A \in$ span$(S)'$ if it satisfies the following Fréchet-like differentiability condition for all $f, g \in$ span$(S)$:

$$\lim_{h \to 0} \frac{1}{h}(A(f + hg) - A(f)) = \nabla A(f)^T g,$$

where $\nabla A(f)^T g$ denotes the value of the linear functional $\nabla A(f)$ at $g$. Note that we adopt the notation $f^T g$ from linear algebra, where it is just the scalar product of the two vectors.

For reference, we shall state the following assumption, which is required in our analysis.

ASSUMPTION 3.1.   Let $A(f)$ be a convex function of $f$ defined on span$(S)$, which satisfies the following conditions:

1. The functional $A$ is differentiable with gradient $\nabla A$.
2. For all $f \in$ span$(S)$ and $g \in S$, the real-valued function $A_{f,g}$ is second-order differentiable (as a function of $h$) and the second derivative satisfies

$$(5) \qquad\qquad A''_{f,g}(0) \le M(\|f\|_1),$$

   where $M(\cdot)$ is a nondecreasing real-valued function.

REMARK 3.1.   A more general form of (5) is $A''_{f,g}(0) \le \ell(g)M(\|f\|_1)$, where $\ell(g)$ is an appropriate scaling factor of $g$. For example, in the examples given below, $\ell(g)$ can be measured by $\sup_x |g(x)|$ or $E_X g(X)^2$. In (5) we assume that functions in $S$ are properly scaled so that $\ell(g) \le 1$. This is for notational convenience only. With more complicated notation techniques developed in this paper can also handle the general case directly without any normalization assumption of the basis functions.

The function $M(\cdot)$ will appear in the convergence analysis in Section 4.1. Although our analysis can handle unbounded $M(\cdot)$, the most interesting boosting examples have bounded $M(\cdot)$ (as we will show shortly). In this case we will also use $M$ to denote a real-valued upper bound of $\sup_a M(a)$.

For statistical estimation problems such as classification and regression with a covariate or predictor variable $X$ and a real response variable $Y$ having a joint distribution, we are interested in the following form of $A(f)$ in (2):

$$(6) \qquad A(f) = \psi(E_{X,Y}\phi(f(X), Y)),$$

where $\phi(\cdot, \cdot)$ is a loss function that is convex in its first argument and $\psi$ is a monotonic increasing auxiliary function which is introduced so that $A(f)$ is convex and $M(\cdot)$ behaves nicely (e.g., bounded). We note that the introduction of $\psi$ is for proving numerical convergence results using our proof techniques, which are needed for proving statistical consistency of boosting with early stopping. However, $\psi$ is not necessary for the actual implementation of the boosting procedure. Clearly the minimizer of (6) that solves (2) does not depend on the choice of $\psi$. Moreover, the behavior of Algorithm 2.1 is not affected by the choice of $\psi$ as long as $\varepsilon_k$ in (3) is appropriately redefined. We may thus always take $\psi(u) = u$, but choosing other auxiliary functions can be convenient for certain problems in our analysis since the resulting formulation has a bounded $M(\cdot)$ function (see the examples given below). We have also used $E_{X,Y}$ to indicate the expectation with respect to the joint distribution of $(X, Y)$.

When not explicitly specified, $E_{X,Y}$ can denote the expectation either with respect to the underlying population or with respect to the empirical samples. This makes no difference as far as our convergence analysis in Section 4.1 is concerned. When it is necessary to distinguish an empirical quantity from its population counterpart, we shall denote the former by a hat above the corresponding quantity. For example, $\hat{E}$ denotes the expectation with respect to the empirical samples, and $\hat{A}$ is the function in (6) with $E_{X,Y}$ replaced by $\hat{E}_{X,Y}$. This distinction will become necessary in the uniform convergence analysis of Section 4.2.

An important application of boosting is binary classification. In this case it is very natural for us to use a set of basis functions that satisfy the conditions

$$(7) \qquad \sup_{g \in S, x} |g(x)| \le 1, \qquad y = \pm 1.$$

For certain loss functions (such as least squares) this condition can be relaxed. In the classification literature $\phi(f, y)$ usually has a form $\phi(fy)$.

Commonly used loss functions are listed in Section A.1. They show that for a typical boosting loss function $\phi$, there exists a constant $M$ such that $\sup_a M(a) \le M$.

3.2. *Assumptions for the statistical convergence analysis.* In classification or regression problems with a covariate or predictor variable $X$ on $R^d$ and a real response variable $Y$, we observe $m$ i.i.d. samples $Z_1^m = \{(X_1, Y_1), \ldots, (X_m, Y_m)\}$ from an unknown underlying distribution $D$. Consider a loss function $\phi(f, y)$ and define $Q(f)$ (true risk) and $\hat{Q}(f)$ (empirical risk) as

$$(8) \quad Q(f) = E_D \phi(f(X), Y), \qquad \hat{Q}(f) = \hat{E}\phi(f(X), Y) = \frac{1}{m} \sum_{i=1}^{m} \phi(f(X_i), Y_i),$$

where $E_D$ is the expectation over the unknown true joint distribution $D$ of $(X, Y)$ (denoted by $E_{X,Y}$ previously); $\hat{E}$ is the empirical expectation based on the sample $Z_1^m$.

Boosting estimators are constructed by applying Algorithm 2.1 with respect to the empirical expectation $\hat{E}$ with a set $S$ of real-valued basis functions $g(x)$. We use $\hat{A}(f)$ to denote the empirical objective function,

$$\hat{A}(f) = \psi(\hat{Q}(f)) = \psi(\hat{E}\phi(f(X), Y)).$$

Similarly, quantities $f_k$, $\alpha_k$ and $g_k$ in Algorithm 2.1 will be replaced by $\hat{f}_k$, $\hat{\alpha}_k$ and $\hat{g}_k$, respectively.

Techniques from modern empirical process theory can be used to analyze the statistical convergence of a boosting estimator with a finite sample. In particular, we use the concept of Rademacher complexity, which is given by the following definition.

DEFINITION 3.2. Let $G = \{g(x, y)\}$ be a set of functions of input $(x, y)$. Let $\{\sigma_i\}_{i=1}^m$ be a sequence of binary random variables such that $\sigma_i = \pm 1$ with probability $1/2$. The (one-sided) sample-dependent *Rademacher complexity* of $G$ is given by

$$R_m(G, Z_1^m) = E_\sigma \sup_{g \in G} \frac{1}{m} \sum_{i=1}^{m} \sigma_i g(X_i, Y_i),$$

and the expected Rademacher complexity of $G$ is denoted by

$$R_m(G) = E_{Z_1^m} R_m(G, Z_1^m).$$

The Rademacher complexity approach for analyzing boosting algorithms first appeared in [21], and it has been used by various people to analyze learning problems, including boosting; for example, see [3, 2, 4, 6, 30]. The analysis using Rademacher complexity as defined above can be applied both to regression and to classification. However, for notational simplicity we focus only on boosting methods for classification, where we impose the following assumption. This assumption is not essential to our analysis, but it simplifies the calculations and some of the final conditions.

ASSUMPTION 3.2. We consider the following form of $\phi$ in (8): $\phi(f, y) = \phi(fy)$ with a convex function $\phi(a): R \to R$ such that $\phi(-a) > \phi(a)$ for all $a > 0$. Moreover, we assume that

(i) Condition (7) holds.
(ii) $S$ in Algorithm 2.1 is closed under negation (i.e., $f \in S \to -f \in S$).
(iii) There exists a finite Lipschitz constant $\gamma_\phi(\beta)$ of $\phi$ in $[-\beta, \beta]$:

$$\forall |f_1|, |f_2| \leq \beta \qquad |\phi(f_1) - \phi(f_2)| \leq \gamma_\phi(\beta)|f_1 - f_2|.$$

The Lipschitz condition of a loss function is usually easy to estimate. For reference, we list $\gamma_\phi$ for loss functions considered in Section A.1:

(a) Logistic regression $\phi(f) = \ln(1 + \exp(-f)): \gamma_\phi(\beta) \leq 1$.
(b) Exponential $\phi(f) = \exp(-f): \gamma_\phi(\beta) \leq \exp(\beta)$.
(c) Least squares $\phi(f) = (f - 1)^2: \gamma_\phi(\beta) \leq 2(\beta + 1)$.
(d) Modified least squares $\phi(f) = \max(1 - f, 0)^2: \gamma_\phi(\beta) \leq 2(\beta + 1)$.
(e) $p$-norm $\phi(f) = |f - 1|^p (p \geq 2): \gamma_\phi(\beta) \leq p(\beta + 1)^{p-1}$.

3.3. *Main statistical results.* We may now state the main statistical results based on the assumptions and definitions given earlier. The following theorem gives conditions for our boosting algorithm so that consistency can be achieved in the large sample limit. The proof is deferred to Section 5.2, with some auxiliary results.

THEOREM 3.1. *Under Assumption 3.2 let $\phi$ be one of the loss functions considered in Section A.1. Assume further that in Algorithm 2.1 we choose quantities $f_0$, $\varepsilon_k$ and $\Lambda_k$ to be independent of the sample $Z_1^m$, such that $\sum_{j=0}^\infty \varepsilon_j < \infty$, and $h_k = \sup \Lambda_k$ satisfies (4).*

*Consider two sequences of sample independent numbers $k_m$ and $\beta_m$ such that $\lim_{m \to \infty} k_m = \infty$ and $\lim_{m \to \infty} \gamma_\phi(\beta_m)\beta_m R_m(S) = 0$. Then as long as we stop Algorithm 2.1 at a step $\hat{k}$ based on $Z_1^m$ such that $\hat{k} \geq k_m$ and $\|\hat{f}_{\hat{k}}\|_1 \leq \beta_m$, we have the consistency result*

$$\lim_{m \to \infty} E_{Z_1^m} Q(\hat{f}_{\hat{k}}) = \inf_{f \in \text{span}(S)} Q(f).$$

REMARK 3.2. The choice of $(k_m, \beta_m)$ in the above theorem should not be void, in the sense that for all samples $Z_1^m$ it should be possible to stop Algorithm 2.1 at a point such that the conditions $\hat{k} \geq k_m$ and $\|\hat{f}_{\hat{k}}\|_1 \leq \beta_m$ are satisfied.

In particular, if $\lim_{m \to \infty} R_m(S) = 0$, then we can always find $k_m \leq k_m'$ such that $k_m \to \infty$ and $\gamma_\phi(\beta_m)\beta_m R_m(S) \to 0$ with $\beta_m = \|f_0\|_1 + \sum_{j=0}^{k_m'} h_j$. This choice of $(k_m, \beta_m)$ is valid as we can stop the algorithm at any $\hat{k} \in [k_m, k_m']$.

Similar to the consistency result, we may further obtain some rate of convergence results. This work does not focus on rate of convergence analysis, and results we obtain are not necessarily tight. Before stating a more general and more complicated result, we first present a version for constant step-size logistic boosting, which is much easier to understand.

THEOREM 3.2. *Consider the logistic regression loss function, with basis $S$ which satisfies $R_m(S) \leq \frac{C_S}{\sqrt{m}}$ for some positive constant $C_S$. For each sample size $m$, consider Algorithm 2.1 with $f_0 = 0$, $\sup_k \Lambda_k = h_0(m) \leq 1/\sqrt{m}$ and $\varepsilon_k \leq h_0(m)^2/2$. Assume that we run boosting for $k(m) = \beta_m/h_0(m)$ steps. Then*

$$E_{Z_1^m} Q(\hat{f}_{\hat{k}}) \leq \inf_{\bar{f} \in \operatorname{span}(S)} \Big[ Q(\bar{f}) + \frac{(2C_S+1)\beta_m}{\sqrt{m}} + \frac{\|\bar{f}\|_1 + 1}{\sqrt{m}} + \frac{\|\bar{f}\|_1}{\|\bar{f}\|_1 + \beta_m} \Big].$$

Note that the condition $R_m(S) \leq C_S/\sqrt{m}$ is satisfied for many basis function classes, such as two-level neural networks and tree basis functions (see Section 4.3). The bound in Theorem 3.2 is independent of $h_0(m)$ [as long as $h_0(m) \leq m^{-1/2}$]. Although this bound is likely to be suboptimal for practice problems, it does give a worst case guarantee for boosting with the greedy optimization aspect taken into consideration. Assume that there exists $\bar{f} \in \operatorname{span}(S)$ such that $Q(\bar{f}) = \inf_{f \in \operatorname{span}} Q(f)$. Then we may choose $\beta_m$ as $\beta_m = O(\|\bar{f}\|_1^{1/2} m^{1/4})$, which gives a convergence rate of $E_{Z_1^m} Q(\hat{f}_{\hat{k}}) \leq Q(\bar{f}) + O(\|\bar{f}\|_1^{1/2} m^{-1/4})$. As the target complexity $\|\bar{f}\|_1$ increases, the convergence becomes slower. An example is provided in Section 6 to illustrate this phenomenon.

We now state the more general result, on which Theorem 3.2 is based (see Section 5.3).

THEOREM 3.3. *Under Assumption 3.2, let $\phi(f) \geq 0$ be a loss function such that $A(f)$ satisfies Assumption 3.1 with the choice $\psi(a) = a$. Given a sample size $m$, we pick a positive nonincreasing sequence $\{h_k\}$ which may depend on $m$. Consider Algorithm 2.1 with $f_0 = 0$, $\sup_k \Lambda_k = h_k$ and $\varepsilon_k \leq h_k^2 M(s_{k+1})/2$, where $s_k = \sum_{i=0}^{k-1} h_i$.*

*Given training data, suppose we run boosting for $\hat{k} = k(m)$ steps, and let $\beta_m = s_{k(m)}$. Then $\forall \bar{f} \in \operatorname{span}(S)$ such that $Q(\bar{f}) \leq Q(0)$*

$$E_{Z_1^m} Q(\hat{f}_{\hat{k}}) \leq Q(\bar{f}) + 2\gamma_\phi(\beta_m)\beta_m R_m(S)$$
$$+ \frac{1}{\sqrt{m}}\phi(-\|\bar{f}\|_1) + \frac{\|\bar{f}\|_1 \phi(0)}{\|\bar{f}\|_1 + \beta_m} + \delta_m(\|\bar{f}\|_1),$$

*where*

$$\delta_m(\|\bar{f}\|_1) = \inf_{1 \le \ell \le k(m)} \left[ \frac{s_\ell + \|\bar{f}\|_1}{\beta_m + \|\bar{f}\|_1} h_0^2 \ell + (k(m) - \ell) h_\ell^2 \right] M(\beta_m + h_{k(m)}).$$

If the target function is $\bar{f}$ which belongs to span($S$), then Theorem 3.3 can be directly interpreted as a rate of convergence result. However, the expression of $\delta_m$ may still be quite complicated. For specific loss function and step-size choices, the bound can be simplified. For example, the result for logistic boosting in Theorem 3.2 follows easily from the theorem (see Section 5.3).

**4. Preparatory results.** As discussed earlier, it is well known by now that boosting can overfit if left to run until convergence. In Section 3.3 we stated our main results that with appropriately chosen stopping rules and under regularity conditions, results of consistency and rates of convergence can be obtained. In this section we begin the proof process of these main results by proving the necessary preparatory results, which are interesting in their own right, especially those on numerical convergence of boosting in Section 4.1.

Suppose that we run Algorithm 2.1 on the sample $Z_1^m$ and stop at step $\hat{k}$. By the triangle inequality and for any $\bar{f} \in \text{span}(S)$, we have

$$
\begin{aligned}
(9) \qquad E_{Z_1^m} Q(\hat{f}_{\hat{k}}) - Q(\bar{f}) &\le E_{Z_1^m} |\hat{Q}(\hat{f}_{\hat{k}}) - Q(\hat{f}_{\hat{k}})| + E_{Z_1^m} |\hat{Q}(\bar{f}) - Q(\bar{f})| \\
&\quad + E_{Z_1^m} [\hat{Q}(\hat{f}_{\hat{k}}) - \hat{Q}(\bar{f})].
\end{aligned}
$$

The middle term is on a fixed $\bar{f}$, and thus it has a rate of convergence $O(1/\sqrt{m})$ by the CLT. To study the consistency and rates of convergence of boosting with early stopping, the work lies in dealing with the first and third terms in (9). The third term is on the empirical performance of the boosting algorithm, and thus a numerical convergence analysis is required and hence proved in Section 4.1. Using modern empirical process theory, in Section 4.2 we upper bound the first term in terms of Rademacher complexity.

We will focus on the loss functions (such as those in Section A.1) which satisfy Assumption 3.1. In particular, we assume that $\psi$ is a monotonic increasing function, so that minimizing $A(f)$ or $\hat{A}(f)$ is equivalent to minimizing $Q(f)$ or $\hat{Q}(f)$. The derivation in Section 4.2 works with $Q(f)$ and $\hat{Q}(f)$ directly, instead of $A(f)$ and $\hat{A}(f)$. The reason is that, unlike our convergence analysis in Section 4.1, the relatively simple sample complexity analysis presented in Section 4.2 does not take advantage of $\psi$.

4.1. *Numerical convergence analysis.* Here we consider the numerical convergence behavior of $f_k$ obtained from the greedy boosting procedure as $k$ increases. For notational simplicity, we state the convergence results in terms of the population boosting algorithm, even though they also hold for the empirical boosting algorithm. The proofs of the two main lemmas are deferred to Section A.2.

In our convergence analysis, we will specify convergence bounds in terms of $\|\bar{f}\|_1$ (where $\bar{f}$ is a reference function) and a sequence of nondecreasing numbers $s_k$ satisfying the following condition: there exist positive numbers $h_k$ such that

$$(10) \qquad |\bar{\alpha}_k| \le h_k \in \Lambda_k \qquad \text{and let } s_k = \|f_0\|_1 + \sum_{i=0}^{k-1} h_i,$$

where $\{\bar{\alpha}_k\}$ are the step-sizes in (3). Note that $h_k$ in (10) can be taken as any number that satisfies the above condition, and it can depend on $\{\bar{\alpha}_k\}$ computed by the boosting algorithm. However, it is often desirable to state a convergence result that does not depend on the actual boosting outputs (i.e., the actual $\bar{\alpha}_k$ computed). For such results we may simply fix $h_k$ by letting $h_k = \sup \Lambda_k$. This gives convergence bounds for the restricted step-size method which we mentioned earlier.

It can be shown (see Section A.2) that even in the worse case, the value $A(f_{k+1}) - A(\bar{f})$ decreases from $A(f_k) - A(\bar{f})$ by a reasonable quantity. Cascading this analysis leads to a numerical rate or speed of convergence for the boosting procedure.

The following lemma contains the one-step convergence bound, which is the key result in our convergence analysis.

LEMMA 4.1. *Assume that $A(f)$ satisfies Assumption* 3.1. *Consider $h_k$ and $s_k$ that satisfy* (10). *Let $\bar{f}$ be an arbitrary reference function in* $\mathrm{span}(S)$, *and define*

$$(11) \qquad \Delta A(f_k) = \max(0, A(f_k) - A(\bar{f})),$$

$$(12) \qquad \bar{\varepsilon}_k = \frac{h_k^2}{2} M(s_{k+1}) + \varepsilon_k.$$

*Then after $k$ steps, the following bound holds for $f_{k+1}$ obtained from Algorithm* 2.1:

$$(13) \qquad \Delta A(f_{k+1}) \le \left(1 - \frac{h_k}{s_k + \|\bar{f}\|_1}\right) \Delta A(f_k) + \bar{\varepsilon}_k.$$

Applying Lemma 4.1 repeatedly, we arrive at a convergence bound for the boosting Algorithm 2.1 as in the following lemma.

LEMMA 4.2.   *Under the assumptions of Lemma* 4.1*, we have*

$$(14) \qquad \Delta A(f_k) \leq \frac{\|f_0\|_1 + \|\bar{f}\|_1}{s_k + \|\bar{f}\|_1} \Delta A(f_0) + \sum_{j=1}^{k} \frac{s_j + \|\bar{f}\|_1}{s_k + \|\bar{f}\|_1} \bar{\varepsilon}_{j-1}.$$

The above lemma gives a quantitative bound on the convergence of $A(f_k)$ to the value $A(\bar{f})$ of an arbitrary reference function $\bar{f} \in \text{span}(S)$. We can see that the numerical convergence speed of $A(f_k)$ to $A(\bar{f})$ depends on $\|\bar{f}\|_1$ and the accumulated or total step-size $s_k$. Specifically, if we choose $\bar{f}$ such that $A(\bar{f}) \leq A(f_0)$, then it follows from the above bound that

$$
\begin{aligned}
(15) \quad A(f_{k+1}) &\leq A(\bar{f}) \left\{ 1 - \frac{s_0 + \|\bar{f}\|_1}{s_{k+1} + \|\bar{f}\|_1} \right\} + \frac{s_0 + \|\bar{f}\|_1}{s_{k+1} + \|\bar{f}\|_1} A(f_0) \\
&\quad + \sum_{j=0}^{k} \frac{s_{j+1} + \|\bar{f}\|_1}{s_{k+1} + \|\bar{f}\|_1} \bar{\varepsilon}_j.
\end{aligned}
$$

Note that the inequality is automatically satisfied when $A(f_{k+1}) \leq A(\bar{f})$.

Clearly, in order to select $\bar{f}$ to optimize the bound on the right-hand side, we need to balance a trade-off: we may select $\bar{f}$ such that $A(\bar{f})$ (and thus the first term) becomes smaller as we increase $\|\bar{f}\|_1$; however, the other two terms will become large when $\|\bar{f}\|_1$ increases. This bound also reveals the dependence of the convergence on the initial value of the algorithm $f_0$: the closer $A(f_0)$ gets to the infimum of $A$, the smaller the bound. To our knowledge, this is the first convergence bound for greedy boosting procedures with quantitative numerical convergence speed information.

Previous analyses, including matching pursuit for least squares [29], Breiman's analysis [9] of the exponential loss, as well as the Bregman divergence bound in [12] and the analysis of gradient boosting in [31], were all limiting results without any information on the numerical speed of convergence. The key conceptual difference here is that we do not compare to the optimal value directly, but instead, to the value of an arbitrary $\bar{f} \in \text{span}(S)$, so that $\|\bar{f}\|_1$ can be used to measure the convergence speed. This approach is also crucial for problems where $A(\cdot)$ can take $-\infty$ as its infimum, for which a direct comparison will clearly fail (e.g., Breiman's exponential loss analysis requires smoothness assumptions to prevent this $-\infty$ infimum value).

A general limiting convergence result follows directly from the above lemma.

THEOREM 4.1.   *Assume that* $\sum_{j=0}^{\infty} \bar{\varepsilon}_j < \infty$ *and* $\sum_{j=0}^{\infty} h_j = \infty$*; then we have the following optimization convergence result for the greedy boosting algorithm* (2.1)*:*

$$\lim_{k \to \infty} A(f_k) = \inf_{f \in \text{span}(S)} A(f).$$

PROOF. The assumptions imply that $\lim_{k\to\infty} s_k = \infty$. We can thus construct a nonnegative integer-valued function $k \to j(k) \leq k$ such that $\lim_{k\to\infty} s_{j(k)}/s_k = 0$ and $\lim_{k\to\infty} s_{j(k)} = \infty$.

From Lemma 4.2 we obtain for any fixed $\bar{f}$,

$$\Delta A(f_k) \leq \frac{\|f_0\|_1 + \|\bar{f}\|_1}{s_k + \|\bar{f}\|_1} \Delta A(f_0) + \sum_{j=1}^{k} \frac{s_j + \|\bar{f}\|_1}{s_k + \|\bar{f}\|_1} \bar{\varepsilon}_{j-1}$$

$$= o(1) + \sum_{j=1}^{j(k)} \frac{s_j + \|\bar{f}\|_1}{s_k + \|\bar{f}\|_1} \bar{\varepsilon}_{j-1} + \sum_{j=j(k)+1}^{k} \frac{s_j + \|\bar{f}\|_1}{s_k + \|\bar{f}\|_1} \bar{\varepsilon}_{j-1}$$

$$\leq o(1) + \frac{s_{j(k)} + \|\bar{f}\|_1}{s_k + \|\bar{f}\|_1} \sum_{j=1}^{j(k)} \bar{\varepsilon}_{j-1} + \sum_{j=j(k)+1}^{k} \bar{\varepsilon}_{j-1} = o(1).$$

Therefore $\lim_{k\to\infty} \max(0, A(f_k) - A(\bar{f})) = 0$. Since our analysis applies to any $\bar{f} \in \text{span}(S)$, we can choose $\bar{f}_j \in \text{span}(S)$ such that $\lim_j A(\bar{f}_j) = \inf_{f \in \text{span}(S)} A(f)$. Now from $\lim_{k\to\infty} \max(0, A(f_k) - A(\bar{f}_j)) = 0$, we obtain the theorem. □

COROLLARY 4.1. *For loss functions such as those in Section* A.1, *we have* $\sup_a M(a) < \infty$. *Therefore as long as there exist* $h_j$ *in* (10) *and* $\varepsilon_j$ *in* (3) *such that* $\sum_{j=0}^{\infty} h_j = \infty$, $\sum_{j=0}^{\infty} h_j^2 < \infty$ *and* $\sum_{j=0}^{\infty} \varepsilon_j < \infty$, *we have the following convergence result for the greedy boosting procedure:*

$$\lim_{k\to\infty} A(f_k) = \inf_{f \in \text{span}(S)} A(f).$$

The above results regarding population minimization automatically apply to the empirical minimization if we assume that the starting point $f_0$, as well as quantities $\varepsilon_k$ and $\Lambda_k$ in (3), are sample-independent, and the restricted step-size case where $h_k = \sup \Lambda_k$ satisfies the condition (4).

The idea of restricting the step-size when we compute $\bar{\alpha}_j$ was advocated by Friedman, who discovered empirically that taking small step-size helps [14]. In our analysis, we can restrict the search region so that Corollary 4.1 is automatically satisfied. Since we believe this is an important case which applies for general loss functions, we shall explicitly state the corresponding convergence result below.

COROLLARY 4.2. *Consider a loss function (e.g., those in Section* A.1*) such that* $\sup_a M(a) < +\infty$. *Pick any sequence of positive numbers* $h_j$ $(j \geq 0)$ *such that* $\sum_{j=0}^{\infty} h_j = \infty$, $\sum_{j=0}^{\infty} h_j^2 < \infty$. *If we choose* $\Lambda_k$ *in Algorithm* 2.1 *such that* $h_k = \sup \Lambda_k$, *and* $\varepsilon_j$ *in* (3) *such that* $\sum_{j=0}^{\infty} \varepsilon_j < \infty$, *then*

$$\lim_{k\to\infty} A(f_k) = \inf_{f \in \text{span}(S)} A(f).$$

Note that the above result requires that the step-size $h_j$ be small ($\sum_{j=0}^{\infty} h_j^2 < \infty$), but also not too small ($\sum_{j=0}^{\infty} h_j = \infty$). As discussed above, the first condition prevents large oscillation. The second condition is needed to ensure that $f_k$ can cover the whole space span($S$).

The above convergence results are limiting results that do not carry any convergence speed information. Although with specific choices of $h_k$ and $s_k$ one may obtain such information from (14), the second term on the right-hand side is typically quite complicated. It is thus useful to state a simple result for a specific choice of $h_k$ and $s_k$, which yields more explicit convergence information.

COROLLARY 4.3. *Assume that $A(f)$ satisfies Assumption* 3.1. *Pick a sequence of nonincreasing positive numbers $h_j$ ($j \geq 0$). Suppose we choose $\Lambda_k$ in Algorithm* 2.1 *such that $h_k = \sup \Lambda_k$, and choose $\varepsilon_k$ in* (3) *such that $\varepsilon_k \leq h_k^2 M(s_{k+1})/2$. If we start Algorithm* 2.1 *with $f_0 = 0$, then*

$$\Delta A(f_k) \leq \frac{\|\bar{f}\|_1}{s_k + \|\bar{f}\|_1} \Delta A(f_0) + \inf_{1 \leq \ell \leq k} \left[ \frac{\ell(s_\ell + \|\bar{f}\|_1)}{s_k + \|\bar{f}\|_1} h_0^2 + (k - \ell) h_\ell^2 \right] M(s_{k+1}).$$

PROOF. Using notation of Lemma 4.1, we have $\bar{\varepsilon}_\ell \leq h_\ell^2 M(s_{k+1})$. Therefore each summand in the second term on the right-hand size of Lemma 4.2 is no more than $h_\ell^2 M(s_{k+1})$ when $j > \ell$ and is no more than $h_0^2 M(s_{k+1})(s_\ell + \|\bar{f}\|_1)/(s_k + \|\bar{f}\|_1)$ when $j \leq \ell$. The desired inequality is now a straightforward consequence of (14). □

Note that similar to the proof of Theorem 4.1, the term $(k - \ell)h_\ell^2$ in Corollary 4.3 can also be replaced by $\sum_{j=\ell+1}^{k} h_j^2$. A special case of Corollary 4.3 is constant step-size ($h_k = h_0$) boosting, which is the original version of restricted step-size boosting considered by Friedman [14]. This method is simple to apply since there is only one step-size parameter to choose. Corollary 4.3 shows that boosting with constant step-size (also referred to as $\varepsilon$-boosting in the literature) converges to the optimal value in the limit of $h_0 \to 0$, as long as we choose the number of iterations $k$ and step-size $h_0$ such that $kh_0 \to \infty$ and $kh_0^2 \to 0$. To the best of our knowledge, this is the only rigorously stated convergence result for the $\varepsilon$-boosting method, which justifies why one needs to use a step-size that is as small as possible.

It is also possible to handle sample-dependent choices of $\Lambda_k$ in Algorithm 2.1, or allow unrestricted step-size ($\Lambda_k = R$) for certain formulations. However, the corresponding analysis becomes much more complicated. According to Friedman [14], the restricted step-size boosting procedure is preferable in practice. Therefore we shall not provide a consistency analysis for unrestricted step-size formulations in this paper; but see Section A.3 for relaxations of the restricted step-size condition.

In addition to the above convergence results for general boosting algorithms, Lemma 4.2 has another very useful consequence regarding the limiting behavior of AdaBoost in the separable classification case. It asserts that the infinitely small step-size version of AdaBoost, in the convergence limit, is an $L_1$ margin maximizer. This result has been observed through a connection between boosting with early stopping and $L_1$ constrained boosting (see [18]). Our analysis gives a direct and rigorous proof. This result is interesting because it shows that AdaBoost shares some similarity (in the limit) with support vector machines (SVMs) whose goal in the separable case is to find maximum margin classifiers; the concept of margin has been popularized by Vapnik [36] who used it to analyze the generalization performance of SVMs. The detailed analysis is provided in Section A.4.

4.2. *Uniform convergence.* There are a number of possible ways to study the uniform convergence of empirical processes. In this section we use a relatively simple approach based on Rademacher complexity. Examples with neural networks and tree-basis (left orthants) functions will be given to illustrate our analysis.

The Rademacher complexity approach for analyzing boosting algorithms appeared first in [21]. Due to its simplicity and elegance, it has been used and generalized by many researchers [2, 3, 4, 6, 30]. The approach used here essentially follows Theorem 1 of [21], but without concentration results.

From Lemma 4.2 we can see that the convergence of the boosting procedure is closely related to $\|\bar{f}\|_1$ and $\|f_k\|_1$. Therefore it is natural for us to measure the learning complexity of Algorithm 2.1 based on the 1-norm of the function family it can approximate at any given step. We shall mention that this analysis is not necessarily the best approach for obtaining tight learning bounds since the boosting procedure may effectively search a much smaller space than the function family measured by the 1-norm $\|f_k\|_1$. However, it is relatively simple, and sufficient for our purpose of providing an early-stopping strategy to give consistency and some rate of convergence results.

Given any $\beta > 0$, we now would like to estimate the rate of uniform convergence,

$$R_m^\beta = E_{Z_1^m} \sup_{\|f\|_1 \le \beta} (Q(f) - \hat{Q}(f)),$$

where $Q$ and $\hat{Q}$ are defined in (8).

The concept of Rademacher complexity used in our analysis is given in Definition 3.2. For simplicity, our analysis also employs Assumption 3.2. As mentioned earlier, the conditions are not essential, but rather they simplify the final results. For example, the condition (7) implies that $\forall f \in \text{span}(S)$, $|f(x)| \le \|f\|_1$. It follows that $\forall \beta \ge \|f\|_1$, $\phi(f, y) \le \phi(-\beta)$. This inequality, although convenient, is certainly not essential.

LEMMA 4.3. *Under Assumption* 3.2,

$$R_m^\beta = E_{Z_1^m} \sup_{\|f\|_1 \le \beta} [E_D \phi(f(X), Y) - \hat{E}\phi(f(X), Y)] \le 2\gamma_\phi(\beta)\beta R_m(S), \tag{16}$$

*where* $\gamma_\phi(\beta)$ *is a Lipschitz constant of* $\phi$ *in* $[-\beta, \beta]$: $\forall |f_1|, |f_2| \le \beta$: $|\phi(f_1) - \phi(f_2)| \le \gamma_\phi(\beta)|f_1 - f_2|$.

PROOF. Using the standard symmetrization argument (e.g., see Lemma 2.3.1 of [35]), we have

$$R_m^\beta = E_{Z_1^m} \sup_{\|f\|_1 \le \beta} [E_D \phi(f(X), Y) - \hat{E}\phi(f(X), Y)]$$

$$\le 2R_m(\{\phi(f(X), Y) : \|f\|_1 \le \beta\}).$$

Now the one-sided Rademacher process comparison result in [32], Theorem 7, which is essentially a slightly refined result (with better constant) of the two-sided version in [24], Theorem 4.12, implies that

$$R_m(\{\phi(f(X), Y) : \|f\|_1 \le \beta\}) \le \gamma_\phi(\beta) R_m(\{f(X) : \|f\|_1 \le \beta\}).$$

Using the simple fact that $g = \sum_i \alpha_i f_i$ ($\sum_i |\alpha_i| = 1$) implies $g \le \max(\sup_i f_i, \sup_i -f_i)$, and that $S$ is closed under negation, it is easy to verify that $R_m(S) = R_m(\{f \in \text{span}(S) : \|f\|_1 \le 1\})$. Therefore

$$R_m(\{f(X) : \|f\|_1 \le \beta\}) = \beta R_m(S).$$

Now by combining the three inequalities, we obtain the lemma. □

4.3. *Estimating Rademacher complexity.* Our uniform convergence result depends on the Rademacher complexity $R_m(S)$. For many function classes, it can be estimated directly. In this section we use a relation between Rademacher complexity and $\ell_2$-covering numbers from [35].

Let $X = \{X_1, \ldots, X_m\}$ be a set of points and let $Q_m$ be the uniform probability measure over these points. We define the $\ell_2(Q_m)$ distance between any two functions $f$ and $g$ as

$$\ell_2(Q_m)(f, g) = \left(\frac{1}{m} \sum_{i=1}^m |f(x_i) - g(x_i)|^2\right)^{1/2}.$$

Let $F$ be a class of functions. The *empirical $\ell_2$-covering number* of $F$, denoted by $N(\varepsilon, F, \ell_2(Q_m))$, is the minimal number of balls $\{g : \ell_2(Q_m)(g, f) \le \varepsilon\}$ of radius $\varepsilon$ needed to cover $F$. The *uniform $\ell_2$ covering number* is given by

$$N_2(\varepsilon, F, m) = \sup_{Q_m} N(\varepsilon, F, \ell_2(Q_m)),$$

where the supremum is over all probability distribution $Q_m$ over samples of size $m$. If $F$ contains 0, then there exists a universal constant $C$ (see Corollary 2.2.8 in [35]) such that

$$R_m(F) \leq \left( \int_0^\infty \sqrt{\log N_2(\varepsilon, F, m)} \, d\varepsilon \right) \frac{C}{\sqrt{m}},$$

where we assume that the integral on the right-hand side is finite. Note that for a function class $F$ with divergent integration value on the right-hand side, the above inequality can be easily modified so that we start the integration from a point $\varepsilon_0 > 0$ instead of 0. However, the dependency of $R_m(F)$ on $m$ can be slower than $1/\sqrt{m}$.

ASSUMPTION 4.1.    $F$ satisfies the condition

$$\sup_m \int_0^\infty \sqrt{\log N_2(\varepsilon, F, m)} \, d\varepsilon < \infty.$$

A function class $F$ that satisfies Assumption 4.1 is also a Donsker class, for which the central limit theorem holds. In statistics and machine learning, one often encounters function classes $F$ with finite VC-dimension, where the following condition holds (see Theorem 2.6.7 of [35]) for some constants $C$ and $V$ independent of $m$: $N_2(\varepsilon, F, m) \leq C(1/\varepsilon)^V$. Clearly a function class with finite VC-dimension satisfies Assumption 4.1.

For simplicity, in this paper we assume that $S$ satisfies Assumption 4.1. It follows that

$$(17) \qquad\qquad R_m(S) \leq R_m(S \cup \{0\}) \leq \frac{C_S}{\sqrt{m}},$$

where $C_S$ is a constant that depends on $S$ only. This is the condition used in Theorem 3.2. We give two examples of basis functions that are often used in practice with boosting.

*Two-level neural networks.*    We consider two-level neural networks in $R^d$, which form the function space span($S$) with $S$ given by

$$S = \{\sigma(w^T x + b) : w \in R^d, b \in R\},$$

where $\sigma(\cdot)$ is a monotone bounded continuous activation function.

It is well known that $S$ has a finite VC-dimension, and thus satisfies Assumption 4.1. In addition, for any compact subset $U \in R^d$, it is also well known that span($S$) is dense in $C(U)$ (see [26]).

*Tree-basis functions.* Tree-basis (left orthant) functions in $R^d$ are given by the indicator function of rectangular regions,

$$S = \{I((-\infty, a_1] \times \cdots \times (-\infty, a_d]) : a_1, \ldots, a_d \in R\}.$$

Similar to two-level neural networks, it is well known that $S$ has a finite VC-dimension, and for any compact set $U \in R^d$, span$(S)$ is dense in $C(U)$.

In addition to rectangular region basis functions, we may also consider a basis $S$ consisting of restricted size classification and regression trees (disjoint unions of constant functions on rectangular regions), where we assume that the number of terminal nodes is no more than a constant $V$. Such a basis set $S$ also has a finite VC-dimension.

## 5. Consistency and rates of convergence with early stopping.
In this section we put together the results in the preparatory Section 4 to prove consistency and some rate of convergence results for Algorithm 2.1 as stated in the main result Section 3.3. For simplicity we consider only restricted step-size boosting with relatively simple strategies for choosing step-sizes. According to Friedman [14], the restricted step-size boosting procedure is preferable in practice. Therefore we shall not provide a consistency analysis for unrestricted step-size formulations in this paper. Discussions on the relaxation of the step-size condition can be found in Section A.3.

### 5.1. *General decomposition.*
Suppose that we run the boosting algorithm and stop at an early stopping point $\hat{k}$. The quantity $\hat{k}$, which is to be specified in Section 5.2, may depend on the empirical sample $Z_1^m$. Suppose also that the stopping point $\hat{k}$ is chosen so that the resulting boosting estimator $\hat{f}_{\hat{k}}$ satisfies

$$\lim_{m \to \infty} E_{Z_1^m} Q(\hat{f}_{\hat{k}}) = \inf_{f \in \text{span}(S)} Q(f), \tag{18}$$

where we use $E_{Z_1^m}$ to denote the expectation with respect to the random sample $Z_1^m$. Since $Q(\hat{f}_{\hat{k}}) \geq \inf_{f \in \text{span}(S)} Q(f)$, we also have

$$\lim_{m \to \infty} E_{Z_1^m} \left| Q(\hat{f}_{\hat{k}}) - \inf_{f \in \text{span}(S)} Q(f) \right| = \lim_{m \to \infty} E_{Z_1^m} Q(\hat{f}_{\hat{k}}) - \inf_{f \in \text{span}(S)} Q(f) = 0.$$

If we further assume there is a unique $f^*$ such that

$$Q(f^*) = \inf_{f \in \text{span}(S)} Q(f),$$

and for any sequence $\{f_m\}$, $Q(f_m) \to Q(f^*)$ implies that $f_m \to f^*$, then since $Q(\hat{f}_{\hat{k}}) \to Q(f^*)$ as $m \to \infty$, it follows that

$$\hat{f}_{\hat{k}} \to f^* \qquad \text{in probability,}$$

which gives the usual consistency of the boosting estimator with an appropriate early stopping if the target function $f$ coincides with $f^*$. This is the case, for example, if the regression function $f(x) = E_D(Y|x)$ with respect to the true distribution $D$ is in $\mathrm{span}(S)$ or can be approximated arbitrarily close by functions in $\mathrm{span}(S)$.

In the following, we derive a general decomposition needed for proving (18) or Theorem 3.1 in Section 3.3. Suppose that Assumption 3.2 holds. Then for all fixed $\bar{f} \in \mathrm{span}(S)$, we have

$$
\begin{aligned}
E_{Z_1^m} |\hat{Q}(\bar{f}) - Q(\bar{f})| &\leq [E_{Z_1^m} |\hat{Q}(\bar{f}) - Q(\bar{f})|^2]^{1/2} \\
&= \left[ \frac{1}{m} E_D |\phi(\bar{f}(X)Y) - Q(\bar{f})|^2 \right]^{1/2} \\
&\leq \left[ \frac{1}{m} E_D \phi(\bar{f}(X)Y)^2 \right]^{1/2} \leq \frac{1}{\sqrt{m}} \phi(-\|\bar{f}\|_1).
\end{aligned}
$$

Assume that we run Algorithm 2.1 on the sample $Z_1^m$ and stop at step $\hat{k}$. If the stopping point $\hat{k}$ satisfies $P(\|\hat{f}_{\hat{k}}\|_1 \leq \beta_m) = 1$ for some sample-independent $\beta_m \geq 0$, then using the uniform convergence estimate in (16), we obtain

$$
\begin{aligned}
(19) \quad & E_{Z_1^m} Q(\hat{f}_{\hat{k}}) - Q(\bar{f}) \\
&= E_{Z_1^m}[Q(\hat{f}_{\hat{k}}) - \hat{Q}(\hat{f}_{\hat{k}})] + E_{Z_1^m}[\hat{Q}(\bar{f}) - Q(\bar{f})] \\
&\quad + E_{Z_1^m}[\hat{Q}(\hat{f}_{\hat{k}}) - \hat{Q}(\bar{f})] \\
&\leq 2\gamma_\phi(\beta_m)\beta_m R_m(S) + \frac{1}{\sqrt{m}} \phi(-\|\bar{f}\|_1) + \sup_{Z_1^m}[\hat{Q}(\hat{f}_{\hat{k}}) - \hat{Q}(\bar{f})].
\end{aligned}
$$

5.2. *Consistency with restricted step-size boosting.* We consider a relatively simple early-stopping strategy for restricted step-size boosting, where we take $h_k = \sup \Lambda_k$ to satisfy (4).

Clearly, in order to prove consistency, we only need to stop at a point such that $\forall \bar{f} \in \mathrm{span}(S)$, all three terms in (19) become nonpositive in the limit $m \to \infty$. By estimating the third term using Lemma 4.2, we obtain the following proof of our main consistency result (Theorem 3.1).

PROOF OF THEOREM 3.1. Obviously the assumptions of the theorem imply that the first two terms of (19) automatically converge to zero. In the following, we only need to show that $\forall \bar{f} \in \mathrm{span}(S) : \sup_{Z_1^m} \max(0, \hat{Q}(\hat{f}_{\hat{k}}) - \hat{Q}(\bar{f})) \to 0$ when $m \to \infty$.

From Section A.1 we know that there exists a distribution-independent number $M > 0$ such that $M(a) < M$ for all underlying distributions. Therefore for all empirical samples $Z_1^m$, Lemma 4.2 implies that

$$\Delta \hat{A}(\hat{f}_{\hat{k}}) \leq \frac{\|f_0\|_1 + \|\bar{f}\|_1}{s_{\hat{k}} + \|\bar{f}\|_1} \Delta \hat{A}(f_0) + \sum_{j=1}^{\hat{k}} \frac{s_j + \|\bar{f}\|_1}{s_{\hat{k}} + \|\bar{f}\|_1} \bar{\varepsilon}_{j-1},$$

where $\Delta \hat{A}(f) = \max(0, \hat{A}(f) - \hat{A}(\bar{f}))$, $s_k = \|f_0\|_1 + \sum_{i=0}^{k-1} h_i$ and $\bar{\varepsilon}_k = \frac{h_k^2}{2} M + \varepsilon_k$. Now using the inequality $\Delta \hat{A}(f_0) \leq \max(\psi(\phi(-\|f_0\|_1)) - \psi(\phi(\|\bar{f}\|_1)), 0) = c(\bar{f})$ and $\hat{k} \geq k_m$, we obtain

$$(20) \quad \sup_{Z_1^m} \Delta \hat{A}(\hat{f}_{\hat{k}}) \leq \sup_{k \geq k_m} \left[ \frac{\|f_0\|_1 + \|\bar{f}\|_1}{s_k + \|\bar{f}\|_1} c(\bar{f}) + \sum_{j=1}^{k} \frac{s_j + \|\bar{f}\|_1}{s_k + \|\bar{f}\|_1} \bar{\varepsilon}_{j-1} \right].$$

Observe that the right-hand side is independent of the sample $Z_1^m$. From the assumptions of the theorem, we have $\sum_{j=0}^{\infty} \bar{\varepsilon}_j < \infty$ and $\lim_{k \to \infty} s_k = \infty$. Now the proof of Theorem 4.1 implies that as $k_m \to \infty$, the right-hand side of (20) converges to zero. Therefore $\lim_{m \to \infty} \sup_{Z_1^m} \Delta \hat{A}(\hat{f}_{\hat{k}}) = 0$. □

The following universal consistency result is a straightforward consequence of Theorem 3.1.

COROLLARY 5.1. *Under the assumptions of Theorem 3.1, for any Borel set $U \subset R^d$, if span$(S)$ is dense in $C(U)$—the set of continuous functions under the uniform-norm topology, then for all Borel measure $D$ on $U \times \{-1, 1\}$,*

$$\lim_{m \to \infty} E_{Z_1^m} Q(\hat{f}_{\hat{k}}) = \inf_{f \in B(U)} Q(f),$$

*where $B(U)$ is the set of Borel measurable functions.*

PROOF. We only need to show $\inf_{f \in \text{span}(S)} Q(f) = \inf_{f \in B(U)} Q(f)$. This follows directly from Theorem 4.1 of [38]. □

For binary classification problems where $y = \pm 1$, given any real-valued function $f$, we predict $y = 1$ if $f(x) \geq 0$ and $y = -1$ if $f(x) < 0$. The classification error is the following 0–1 loss function:

$$\ell(f(x), y) = I[yf(x) \leq 0],$$

where $I[E]$ is the indicator function of the event $E$, and the expected loss is

$$(21) \quad L(f) = E_D \ell(f(X), Y).$$

The goal of classification is to find a predictor $f$ to minimize (21). Using the notation $\eta(x) = P(Y = 1 | X = x)$, it is well known that $L^*$, the minimum of $L(f)$, can be achieved by setting $f(x) = 2\eta(x) - 1$. Let $D$ be a Borel measure defined on $U \times \{-1, 1\}$; it is known (e.g., see [38]) that if $Q(f) \to \inf_{f \in B(U)} Q(f)$, then $L(f) \to L^*$. We thus have the following consistency result for binary-classification problems.

COROLLARY 5.2. *Under the assumptions of Corollary 5.1, we have*

$$\lim_{m \to \infty} E_{Z_1^m} L(\hat{f}_{\hat{k}}) = L^*.$$

The stopping criterion given in Theorem 3.1 depends on $R_m(S)$. For $S$ that satisfies Assumption 4.1, this can be estimated from (17). The condition $\gamma_\phi(\beta_m)\beta_m R_m(S) \to 0$ in Theorem 3.1 becomes $\gamma_\phi(\beta_m)\beta_m = o(\sqrt{m})$. Using the bounds for $\gamma_\phi(\cdot)$ in Section 4.2, we obtain the following condition.

ASSUMPTION 5.1. The sequence $\beta_m$ satisfies:

(i) Logistic regression $\phi(f) = \ln(1 + \exp(-f))$: $\beta_m = o(m^{1/2})$.
(ii) Exponential $\phi(f) = \exp(-f)$: $\beta_m = o(\log m)$.
(iii) Least squares $\phi(f) = (f - 1)^2$: $\beta_m = o(m^{1/4})$.
(iv) Modified least squares $\phi(f) = \max(0, 1 - f)^2$: $\beta_m = o(m^{1/4})$.
(v) $p$-norm $\phi(f) = |f - 1|^p (p \geq 2)$: $\beta_m = o(m^{1/2p})$.

We can summarize the above discussion in the following theorem, which applies to boosted VC-classes such as boosted trees and two-level neural networks.

THEOREM 5.1. *Under Assumption 3.2, let $\phi$ be one of the loss functions considered in Section A.1. Assume further that in Algorithm 2.1 we choose the quantities $f_0$, $\varepsilon_k$ and $\Lambda_k$ to be independent of the sample $Z_1^m$, such that $\sum_{j=0}^{\infty} \varepsilon_j < \infty$, and $h_k = \sup \Lambda_k$ satisfies (4).*

*Suppose $S$ satisfies Assumption 4.1 and we choose sample-independent $k_m \to \infty$, such that $\beta_m = \|f_0\|_1 + \sum_{j=0}^{k_m} h_j$ satisfies Assumption 5.1. If we stop Algorithm 2.1 at step $k_m$, then $\|\hat{f}_{k_m}\|_1 \leq \beta_m$ and the following consistency result holds:*

$$\lim_{m \to \infty} E_{Z_1^m} Q(\hat{f}_{k_m}) = \inf_{f \in \mathrm{span}(S)} Q(f).$$

*Moreover, if $\mathrm{span}(S)$ is dense in $C(U)$ for a Borel set $U \subset R^d$, then for all Borel measures $D$ on $U \times \{-1, 1\}$, we have*

$$\lim_{m \to \infty} E_{Z_1^m} Q(\hat{f}_{k_m}) = \inf_{f \in B(U)} Q(f), \qquad \lim_{m \to \infty} E_{Z_1^m} L(\hat{f}_{k_m}) = L^*.$$

Note that in the above theorem the stopping criterion $k_m$ is sample-independent. However, similar to Theorem 3.1, we may allow other sample-dependent $\hat{k}$ such that $\|f_{\hat{k}}\|_1$ stays within the $\beta_m$ bound. One may be tempted to interpret the rates of $\beta_m$. However, since different loss functions approximate the underlying distribution in different ways, it is not clear that one can rigorously compare them. Moreover, our analysis is likely to be loose.

5.3. *Some bounds on the rate of convergence.*    In addition to consistency, it is also useful to study statistical rates of convergence of the greedy boosting method with certain target function classes. Since our analysis is based on the 1-norm of the target function, the natural function classes we may consider are those that can be approximated well using a function in span($S$) with small 1-norm.

We would like to emphasize that rate results, that have been stated in Theorems 3.2 and 3.3 and are to be proved here, are not necessarily optimal. There are several reasons for this. First, we relate the numerical behavior of boosting to 1-norm regularization. In reality, this may not always be the best way to analyze boosting since boosting can be studied using other complexity measures such as sparsity (e.g., see [22] for some other complexity measures). Second, even with the 1-norm regularization complexity measure, the numerical convergence analysis in Section 4.1 may not be tight. This again will adversely affect our final bounds. Third, our uniform convergence analysis, based on the relatively simple Rademacher complexity, is not necessarily tight. For some problems there are more sophisticated methods which improve upon our approach here (e.g., see [[2, 3, 4, 5, 6], [22, 30]]).

A related point is that bounds we are interested in here are a priori convergence bounds that are data-independent. In recent years, there has been much interest in developing data-dependent bounds which are tighter (see references mentioned above). For example, in our case we may allow $\beta$ in (16) to depend on the observed data (rather than simply setting it to be a value based only on the sample size). This approach, which can tighten the final bounds based on observation, is a quite significant recent theoretical advance. However, as mentioned above, there are other aspects of our analysis that can be loose. Moreover, we are mainly interested in worst case scenario upper bounds on the convergence behavior of boosting without looking at the data. Therefore we shall not develop data-dependent bounds here.

The statistical convergence behavior of the boosting algorithm relies on its numerical convergence behavior, which can be estimated using (14). Combined with statistical convergence analysis, we can easily obtain our main rate of convergence result in Theorem 3.3.

PROOF OF THEOREM 3.3.   From (19) we obtain

$$E_{Z_1^m} Q(\hat{f}_{\hat{k}}) \leq Q(\bar{f}) + 2\gamma_\phi(\beta_m)\beta_m R_m(S) + \frac{1}{\sqrt{m}}\phi(-\|\bar{f}\|_1) + \sup_{Z_1^m}[\hat{Q}(\hat{f}_{\hat{k}}) - \hat{Q}(\bar{f})].$$

Now we simply apply Corollary 4.3 to bound the last term. This leads to the desired bound.   □

The result for logistic regression in Theorem 3.2 follows easily from Theorem 3.3.

PROOF OF THEOREM 3.2.   Consider logistic regression loss and constant step-size boosting, where $h_k = h_0(m)$. Note that for logistic regression we have $\gamma_\phi(\beta) \leq 1$, $M(a) \leq 1$, $\phi(-\|\bar{f}\|_1) \leq 1 + \|\bar{f}\|_1$ and $\phi(0) \leq 1$. Using these estimates, we obtain from Theorem 3.3,

$$E_{Z_1^m} Q(\hat{f}_{\hat{k}}) \leq Q(\bar{f}) + 2\beta_m R_m(S) + \frac{\|\bar{f}\|_1 + 1}{\sqrt{m}} + \frac{\|\bar{f}\|_1}{\|\bar{f}\|_1 + \beta_m} + \beta_m h_0(m).$$

Using the estimate of $R_m(S)$ in (17), and letting $h_0(m) \leq 1/\sqrt{m}$, we obtain

$$E_{Z_1^m} Q(\hat{f}_{\hat{k}}) \leq Q(\bar{f}) + \frac{(2C_S + 1)\beta_m}{\sqrt{m}} + \frac{\|\bar{f}\|_1 + 1}{\sqrt{m}} + \frac{\|\bar{f}\|_1}{\|\bar{f}\|_1 + \beta_m}.$$

This leads to the claim.   □

**6. Experiments.**   The purpose of this section is not to reproduce the large number of already existing empirical studies on boosting. Although this paper is theoretical in nature, it is still useful to empirically examine various implications of our analysis, so that we can verify they have observable consequences. For this reason our experiments focus mainly on aspects of boosting with early stopping which have not been addressed in previous studies.

Specifically, we are interested in testing consistency and various issues of boosting with early stopping based on our theoretical analysis. As pointed out in [28], experimentally testing consistency is a very challenging task. Therefore, in this section we have to rely on relatively simple synthetic data, for which we can precisely control the problem and the associated Bayes risk. Such an experimental setup serves the purpose of illustrating main insights revealed by our theoretical analyses.

6.1. *Experimental setup.*   In order to fully control the data generation mechanism, we shall use simple one-dimensional examples. A similar experimental setup was also used in [23] to study various theoretical aspects of voting classification methods.

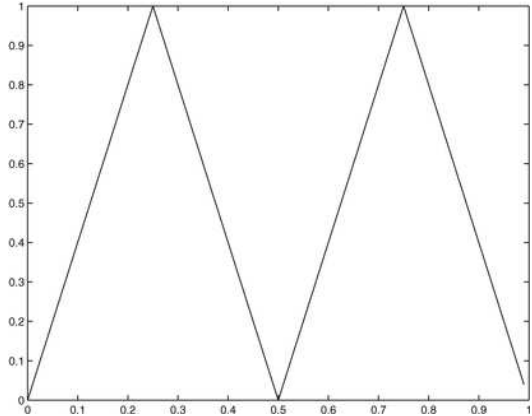

Fig. 1. *Target conditional probability for $d = 2$.*

Our goal is to predict $Y \in \{\pm 1\}$ based on $X \in [0, 1]$. Throughout the experiments, $X$ is uniformly distributed in $[0, 1]$. We consider the target conditional probability of the form $P(Y = 1|X) = 2\{dX\}I(\{dX\} \leq 0.5) + 2(1 - \{dX\}) \times I(\{dX\} > 0.5)$, where $d \geq 1$ is an integer which controls the complexity of the target function, and $I$ denotes the set indicator function. We have also used the notation $\{z\} = z - \lfloor z \rfloor$ to denote the decimal part of a real number $z$, with the standard notation of $\lfloor z \rfloor$ for the integer part of $z$. The Bayes error rate of our model is always 0.25.

Graphically, the target conditional probability contains $d$ triangles. Figure 1 plots such a target for $d = 2$.

We use one-dimensional stumps of the form $I([0, a])$ as our basis functions, where $a$ is a parameter in $[0, 1]$. They form a complete basis since each interval indicator function $I((a, b])$ can be expressed as $I([0, b]) - I([0, a])$.

There have been a number of experimental studies on the impact of using different convex loss functions (e.g., see [14, 27, 28, 39]). Although our theoretical analysis applies to general loss functions, it is not refined enough to suggest that any one particular loss is better than another. For this reason, our experimental study will not include a comprehensive comparison of different loss functions. This task is better left to dedicated empirical studies (such as some of those mentioned above).

We will only focus on consequences of our analysis which have not been well studied empirically. These include various issues related to early stopping and their impact on the performance of boosting. For this purpose, throughout the experiments we shall only use the least-squares loss function. In fact, it is known that this loss function works quite well for many classification problems (see, e.g., [11, 27]) and has been widely applied to many pattern-recognition applications. Its simplicity also makes it attractive.

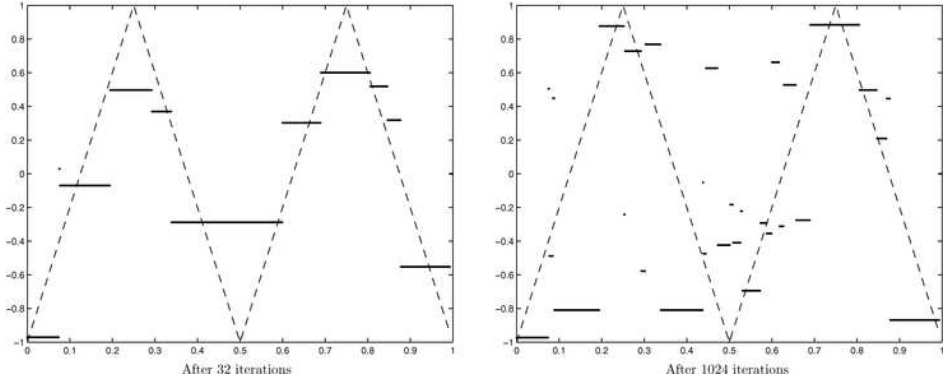

Fig. 2. *Graphs of boosting estimators after $k = 32$ and $1024$ iterations.*

For the least-squares loss, the target function which the boosting procedure tries to estimate is $f_*(x) = 2P(Y = 1|X = x) - 1$. In our experiments, unless otherwise noted, we use boosting with restricted step-size, where at each iteration we limit the step-size to be no larger than $h_i = (i + 1)^{-2/3}$. This choice satisfies our numerical convergence requirement, where we need the conditions $\sum_i h_i = \infty$ and $\sum_i h_i^2 < \infty$. Therefore it also satisfies the consistency requirement in Theorem 3.1.

6.2. *Early stopping and overfitting.* Although it is known that boosting forever can overfit (e.g., see [16, 19]), it is natural to begin our experiments by graphically showing the effect of early-stopping on the predictive performance of boosting.

We shall use the target conditional probability described earlier with complexity $d = 2$, and training sample-size of 100. Figure 2 plots the graphs of estimators obtained after $k = 32$ and 1024 boosting iterations. The dotted lines on the background show the true target function $f_*(x) = 2P(Y = 1|X = x)$. We can see that after 32 iterations, the boosting estimator, although not perfect, roughly has the same shape as that of the true target function. However, after 1024 iterations, the graph appears quite random, implying that the boosting estimator starts to overfit the data.

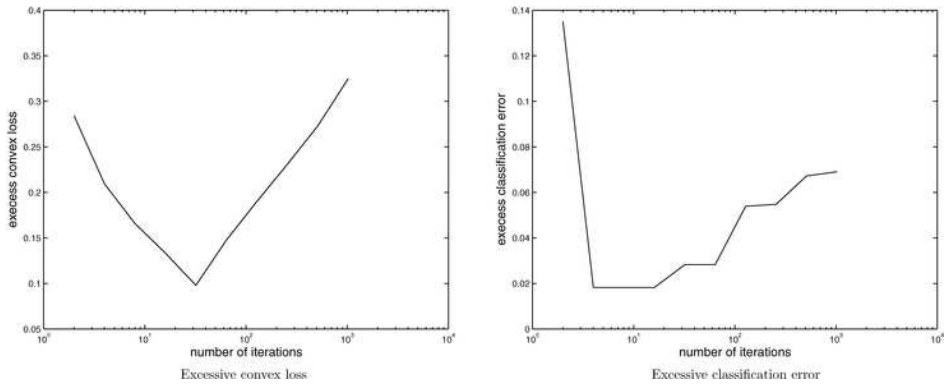

Fig. 3. *Predictive performance of boosting as a function of boosting iterations.*

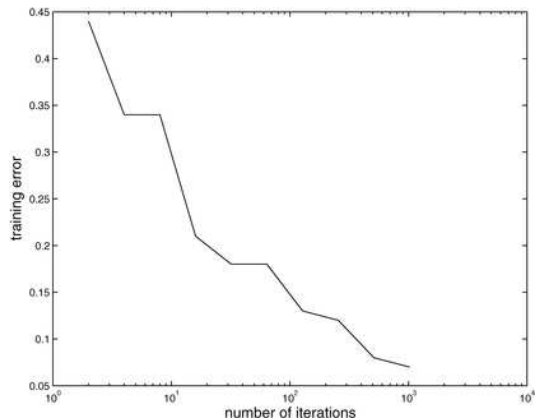

Fig. 4. *Training error.*

Figure 3 shows the predictive performance of boosting versus the number of iterations. The need for early stopping is quite apparent in this example. The *excessive classification error* quantity is defined as the true classification error of the estimator minus the Bayes error (which is 0.25 in our case). Similarly, the *excessive convex loss* quantity is defined as the true least-squares loss of the estimator minus the optimal least-squares loss of the target function $f_*(x)$. Both excessive classification error and convex loss are evaluated through numerical integration for a given decision rule. Moreover, as we can see from Figure 4, the training error continues to decrease as the number of boosting iterations increases, which eventually leads to overfitting of the training data.

6.3. *Early stopping and total step-size.* Since our theoretical analysis favors restricted step-size, a relevant question is what step-size we should choose. We are not the first authors to look into this issue. For example, Friedman and his co-authors suggested using small steps [14, 15]. In fact, they argued that the smaller the step-size, the better. They performed a number of empirical studies to support this claim. Therefore we shall not reinvestigate this issue here. Instead, we focus on a closely related implication of our analysis, which will be useful for the purpose of reporting experimental results in later sections.

Let $\bar{\alpha}_i$ be the step-size taken by the boosting algorithm at the $i$th iteration. Our analysis characterizes the convergence behavior of boosting after

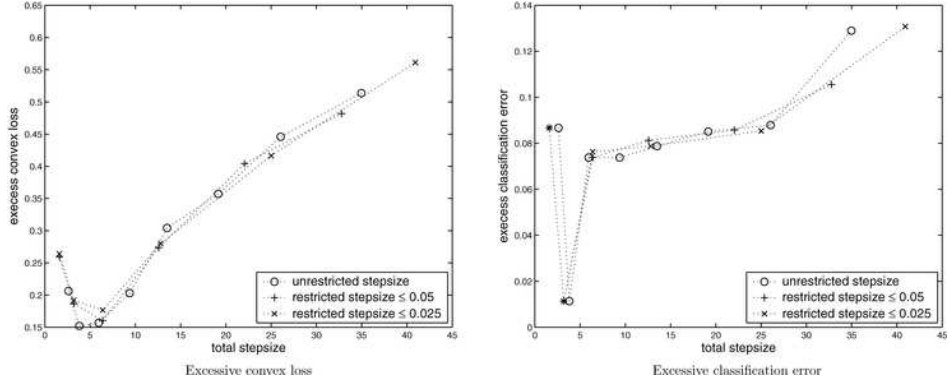

F<small>IG</small>. 5.   *Predictive performance of boosting as a function of total step-size.*

the $k$th step, not by the number of iterations $k$ itself, but rather by the quantity $s_k = \sum_{i \le k} h_i$ in (10), as long as $\bar{\alpha}_i \le h_i \in \Lambda_i$. Although our theorems are stated with the quantity $\sum_{i \le k} h_i$, instead of $\sum_{i \le k} \bar{\alpha}_i$, it does suggest that in order to compare the behavior of boosting under different configurations, it is more natural to use the quantity $\sum_{i \le k} \bar{\alpha}_i$ (which we shall call *total step-size* throughout later experiments) as a measure of stopping point rather than the actual number of boosting iterations. This concept of total step-size also appeared in [18, 17].

Figure 5 shows the predictive performance of boosting versus the total step-size. We use 100 training examples, with the target conditional probability of complexity $d = 3$. The unrestricted step-size method uses exact optimization. Note that for least-squares loss, as explained in Section A.3, the resulting step-sizes will still satisfy our consistency condition $\sum_{i \le k} \bar{\alpha}_i^2 < \infty$. The restricted step-size scheme with step-size $\le h$ employs a constant step-size restriction of $\hat{\alpha}_i \le h$. This experiment shows that the behavior of these different boosting methods is quite similar when we measure the performance not by the number of boosting iterations, but instead by the total step-size. This observation justifies our theoretical analysis, which uses quantities closely related to the total step-size to characterize the convergence behavior of boosting methods. Based on this result, in the next few experiments we shall use the total step-size (instead of the number of boosting iterations) to compare boosting methods under different configurations.

6.4. *The effect of sample-size on early stopping.*   An interesting issue for boosting with early stopping is how its predictive behavior changes when the number of samples increases. Although our analysis does not offer a quantitative characterization, it implies that we should stop later (and the allowable stopping range becomes wider) when sample size increases. This essentially suggests that the optimal stopping point in the boosting predictive performance curve will increase as the sample size increases, and the curve itself becomes flatter. It follows that when the sample size is relatively large, we should run boosting algorithms for a longer time, and it is less necessary to do aggressive early stopping.

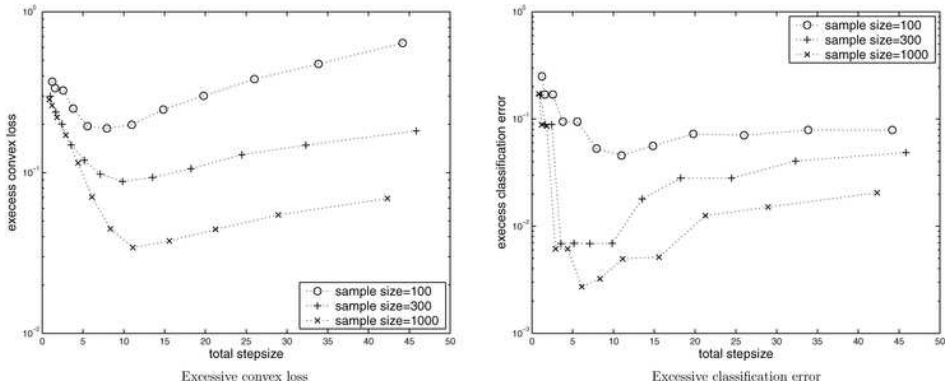

Fig. 6. *Predictive performance of boosting at different sample sizes.*

The above qualitative characterization of the boosting predictive curve also has important practical consequences. We believe this may be one reason why in many practical problems it is very difficult for boosting to overfit, and practitioners often observe that the performance of boosting keeps improving as the number of boosting iterations increases.

Figure 6 shows the effect of sample size on the behavior of the boosting method. Since our theoretical analysis applies directly to the convergence of the convex loss (the convergence of classification error follows implicitly as a consequence of convex loss convergence), the phenomenon described above is more apparent for excessive convex loss curves. The effect on classification error is less obvious, which suggests there is a discrepancy between classification error performance and convex loss minimization performance.

6.5. *Early stopping and consistency.* In this experiment we demonstrate that as sample size increases, boosting with early stopping leads to a consistent estimator with its error rate approaching the optimal Bayes error. Clearly, it is not possible to prove consistency experimentally, which requires running a sample size of $\infty$. We can only use a finite number of samples to demonstrate a clear trend that the predictive performance of boosting with early stopping converges to the Bayes error when the sample size increases. Another main focus of this experiment is to compare the performance of different early stopping strategies.

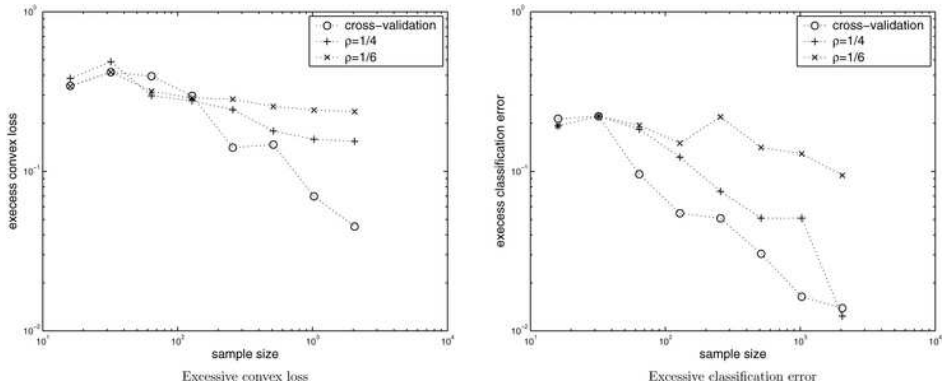

Fig. 7.   *Consistency and early stopping.*

Theoretical results in this paper suggest that for least squares loss, we can achieve consistency as long as we stop at total step-size approximately $m^\rho$ with $\rho < 1/4$, where $m$ is the sample size. We call such an early stopping strategy the $\rho$-*strategy*. Since our theoretical estimate is conservative, we examine the $\rho$-strategy both for $\rho = 1/6$ and for $\rho = 1/4$. Instead of the theoretically motivated (and suboptimal) $\rho$-strategy, in practice one can use cross validation to determine the stopping point. We use a sample size of one-third the training data to estimate the optimal stopping total step-size which minimizes the classification error on the validation set, and then use the training data to compute a boosting estimator which stops at this cross-validation-determined total step-size. This strategy is referred to as the *cross validation* strategy. Figure 7 compares the three early stopping strategies mentioned above. It may not be very surprising to see that the cross-validation-based method is more reliable. The $\rho$-strategies, although they perform less well, also demonstrate a trend of convergence to consistency. We have also noticed that the cross validation scheme stops later than the $\rho$-strategies, implying that our theoretical results impose more restrictive conditions than necessary.

It is also interesting to see how well cross validation finds the optimal stopping point. In Figure 8 we compare the cross validation strategy with two oracle strategies which are not implementable: one selects the optimal stopping point which minimizes the true classification error (which we refer to as *optimal error*), and the other selects the optimal stopping point which minimizes the true convex loss (which we refer to as *optimal convex risk*). These two methods can be regarded as ideal theoretical stopping points for boosting methods. The experiment shows that cross validation performs quite well at large sample sizes.

In the log coordinate space, the convergence curve of boosting with the cross validation stopping criterion is approximately a straight line, which implies that the excess errors decrease as a power of the sample size. By extrapolating this finding, it is reasonable for us to believe that boosting with early stopping converges to the Bayes error in the limit, which verifies the consistency. The two $\rho$ stopping rules, even though showing much slower linear convergence trend, also lead to consistency.

6.6. *The effect of target function complexity on early stopping.* Although we know that boosting with an appropriate early stopping strategy leads to a consistent estimator in the large sample limit, the rate of convergence depends on the complexity of the target function (see Section 5.3). In our analysis the complexity can be measured by the 1-norm of the target function. For target functions considered here, it is not very difficult to show that in order to approximate to an accuracy within $\varepsilon$, it is only necessary to use a combination of our decision stumps with the 1-norm $Cd/\varepsilon$. In this formula $C$ is a constant and $d$ is the complexity of the target function.

Our analysis suggests that the convergence behavior of boosting with early stopping depends on how easy it is to approximate the target function using a combination of basis functions with small 1-norm. A target with $d = u$ is $u$-times as difficult to approximate as a target with $d = 1$. Therefore the optimal stopping point, measured by the total step-size, should accordingly increase as $d$ increases. Moreover, the predictive performance becomes worse. Figure 9 illustrates this phenomenon with $d = 1, 3, 5$ at the sample size of 300. Notice again that since our analysis applies to the convex risk, this phenomenon is much more apparent for the excessive convex loss performance than the excessive classification error performance. Clearly this again shows that although by minimizing a convex loss we indirectly minimize the classification error, these two quantities do not behave identically.

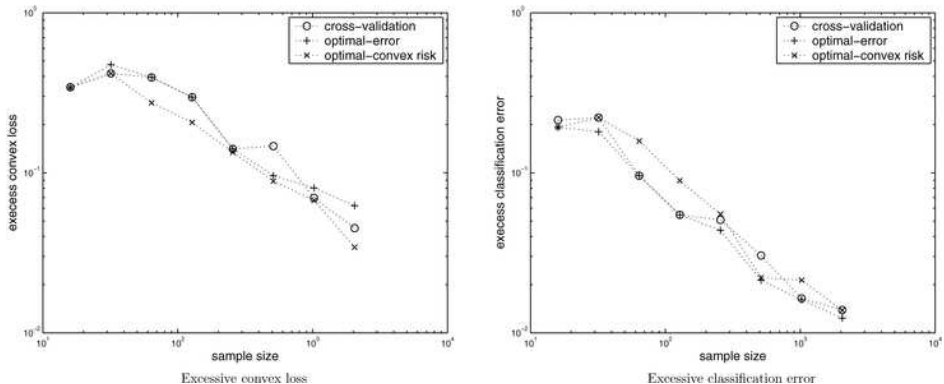

FIG. 8. *Consistency and early stopping.*

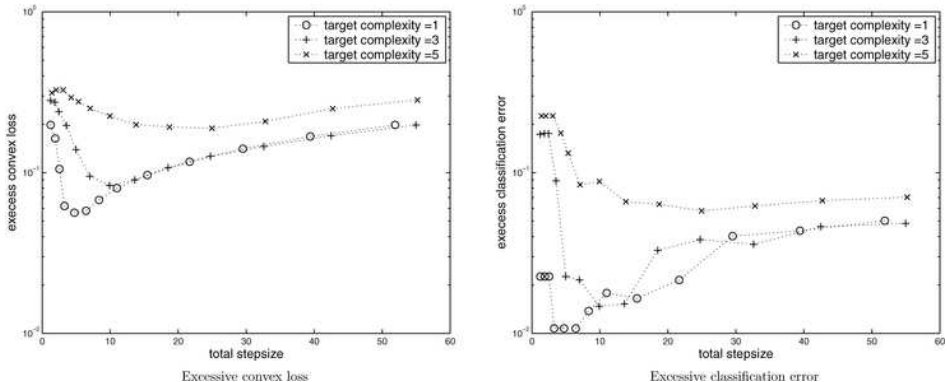

FIG. 9.   *The effect of target function complexity.*

**7. Conclusion.**   In this paper we have studied a general version of the boosting procedure given in Algorithm 2.1. The numerical convergence behavior of this algorithm has been studied using the so-called averaging technique, which was previously used to analyze greedy algorithms for optimization problems defined in the convex hull of a set of basis functions. We have derived an estimate of the numerical convergence speed and established conditions that ensure the convergence of Algorithm 2.1. Our results generalize those in previous studies, such as the matching pursuit analysis in [29] and the convergence analysis of AdaBoost by Breiman [9].

Furthermore, we have studied the learning complexity of boosting algorithms based on the Rademacher complexity of the basis functions. Together with the numerical convergence analysis, we have established a general early stopping criterion for greedy boosting procedures for various loss functions that guarantees the consistency of the obtained estimator in the large sample limit. For specific choices of step-sizes and sample-independent stopping criteria, we have also been able to establish bounds on the statistical rate of convergence. We would like to mention that the learning complexity analysis given in this paper is rather crude. Consequently, the required conditions in our consistency strategy may be more restrictive than one actually needs.

A number of experiments were presented to study various aspects of boosting with early stopping. We specifically focused on issues that have not been covered by previous studies. These experiments show that various quantities and concepts revealed by our theoretical analysis lead to observable consequences. This suggests that our theory can lead to useful insights into practical applications of boosting algorithms.

## APPENDIX

**A.1. Loss function examples.**   We list commonly used loss functions that satisfy Assumption 3.1. They show that for a typical boosting loss function $\phi$, there exists a constant $M$ such that $\sup_a M(a) \leq M$. All loss functions considered are convex.

A.1.1. *Logistic regression.* This is a traditional loss function used in statistics, which is given by (in natural log form here)

$$\phi(f, y) = \ln(1 + \exp(-fy)), \qquad \psi(u) = u.$$

We assume that the basis functions satisfy the condition

$$\sup_{g \in S, x} |g(x)| \leq 1, \qquad y = \pm 1.$$

It can be verified that $A(f)$ is convex differentiable. We also have

$$A''_{f,g}(0) = E_{X,Y} \frac{g(X)^2 Y^2}{(1 + \exp(f(X)Y))(1 + \exp(-f(X)Y))} \leq \frac{1}{4}.$$

A.1.2. *Exponential loss.* This loss function is used in the AdaBoost algorithm, which is the original boosting procedure for classification problems. It is given by

$$\phi(f, y) = \exp(-fy), \qquad \psi(u) = \ln u.$$

Again we assume that the basis functions satisfy the condition

$$\sup_{g \in S, x} |g(x)| \leq 1, \qquad y = \pm 1.$$

In this case it is also not difficult to verify that $A(f)$ is convex differentiable. Hence we also have

$$A''_{f,g}(0) = \frac{E_{X,Y} g(X)^2 Y^2 \exp(-f(X)Y)}{E_{X,Y} \exp(-f(X)Y)} - \frac{[E_{X,Y} g(X)Y \exp(-f(X)Y)]^2}{[E_{X,Y} \exp(-f(X)Y)]^2} \leq 1.$$

A.1.3. *Least squares.* The least squares formulation has been widely studied in regression, but can also be applied to classification problems [10, 11, 14, 30]. A greedy boosting-like procedure for least squares was first proposed in the signal processing community, where it was called *matching pursuit* [29]. The loss function is given by

$$\phi(f, y) = \tfrac{1}{2}(f - y)^2, \qquad \psi(u) = u.$$

We impose the following weaker condition on the basis functions:

$$\sup_{g \in S} E_X g(X)^2 \leq 1, \qquad E_Y Y^2 < \infty.$$

It is clear that $A(f)$ is convex differentiable, and the second derivative is bounded as

$$A''_{f,g}(0) = E_X g(X)^2 \leq 1.$$

A.1.4. *Modified least squares.* For classification problems we may consider the following modified version of the least squares loss, which has a better approximation property [38]:

$$\phi(f,y) = \tfrac{1}{2}\max(1 - fy, 0)^2, \qquad \psi(u) = u.$$

Since this loss is for classification problems, we impose the condition

$$\sup_{g \in S} E_X g(X)^2 \le 1, \qquad y = \pm 1.$$

It is clear that $A(f)$ is convex differentiable, and we have the following bound for the second derivative:

$$A''_{f,g}(0) \le E_X g(X)^2 \le 1.$$

A.1.5. *p-norm boosting.* $p$-norm loss can be interesting both for regression and for classification. In this paper we will only consider the case with $p \ge 2$,

$$\phi(f,y) = |f - y|^p, \qquad \psi(u) = \frac{1}{2(p-1)} u^{2/p}.$$

We impose the condition

$$\sup_{g \in S} E_X |g(X)|^p \le 1, \qquad E_Y |Y|^p < \infty.$$

Now let $u = E_{X,Y}|f(X) + hg(X) - Y|^p$; we have

$$A'_{f,g}(h) = \frac{1}{p-1} u^{(2-p)/p} E_{X,Y} g(X) \operatorname{sign}(f(X) + hg(X) - Y)$$
$$\times |f(X) + hg(X) - Y|^{p-1}.$$

Therefore the second derivative can be bounded as

$$A''_{f,g}(h) = u^{(2-p)/p} E_{X,Y} g(X)^2 |f(X) + hg(X) - Y|^{p-2}$$
$$- \frac{p-2}{p-1} u^{(2-2p)/p}$$
$$\times [E_{X,Y} g(X) \operatorname{sign}(f(X) + hg(X) - Y)|f(X) + hg(X) - Y|^{p-1}]^2$$
$$\le u^{(2-p)/p} E_{X,Y} g(X)^2 |f(X) + hg(X) - Y|^{p-2}$$
$$\le u^{(2-p)/p} E_{X,Y}^{2/p} |g(X)|^p \; E_{X,Y}^{(p-2)/p} |f(X) + hg(X) - Y|^p$$
$$= E_{X,Y}^{2/p} |g(X)|^p \le 1,$$

where the second inequality follows from Hölder's inequality with the duality pair $(p/2, p/(p-2))$.

REMARK A.1. Similar to the least squares case, we can define the modified $p$-norm loss for classification problems. Although the case $p \in (1, 2)$ can be handled by the proof techniques used in this paper, it requires a modified analysis since in this case the corresponding loss function is not second-order differentiable at zero. See related discussions in [37]. Note that the hinge loss used in support vector machines cannot be handled directly with our current technique since its first-order derivative is discontinuous. However, one may approximate the hinge loss with a continuously differentiable function, which can then be analyzed.

**A.2. Numerical convergence proofs.** This section contains two proofs for the numerical convergence analysis section (Section 4.1).

PROOF OF THE ONE-STEP ANALYSIS OR LEMMA 4.1. Given an arbitrary fixed reference function $\bar{f} \in \mathrm{span}(S)$ with the representation

$$\bar{f} = \sum_j \bar{w}^j \bar{f}_j, \qquad \bar{f}_j \in S, \tag{22}$$

we would like to compare $A(f_k)$ to $A(\bar{f})$. Since $\bar{f}$ is arbitrary, we use such a comparison to obtain a bound on the numerical convergence rate.

Given any finite subset $S' \subset S$ such that $S' \supset \{\bar{f}_j\}$, we can represent $\bar{f}$ minimally as

$$\bar{f} = \sum_{g \in S'} \bar{w}_{S'}^g g,$$

where $\bar{w}_{S'}^g = \bar{w}^j$ when $g = \bar{f}_j$ for some $j$, and $\bar{w}_{S'}^g = 0$ when $g \notin \{\bar{f}_j\}$. A quantity that will appear in our analysis is $\|\bar{w}_{S'}\|_1 = \sum_{g \in S'} |\bar{w}_{S'}^g|$. Since $\|\bar{w}_{S'}\|_1 = \|\bar{w}\|_1$, without any confusion, we will still denote $\bar{w}_{S'}$ by $\bar{w}$ with the convention that $\bar{w}^g = 0$ for all $g \notin \{\bar{f}_j\}$.

Given this reference function $\bar{f}$, let us consider a representation of $f_k$ as a linear combination of a finite number of functions $S_k \subset S$, where $S_k \supset \{\bar{f}_j\}$ is to be chosen later. That is, with $g$ indexing an arbitrary function in $S_k$, we expand $f_k$ in terms of $f_k^g$'s which are members of $S_k$ with coefficients $\beta_k^g$:

$$f_k = \sum_{g \in S_k} \beta_k^g f_k^g. \tag{23}$$

With this representation, we define

$$\Delta W_k = \|\bar{w} - \beta_k\|_1 = \sum_{g \in S_k} |\bar{w}^g - \beta_k^g|.$$

Recall that in the statement of the lemma, the convergence bounds are in terms of $\|\bar{w}\|_1$ and a sequence of nondecreasing numbers $s_k$, which satisfy

the condition

$$s_k = \|f_0\|_1 + \sum_{i=0}^{k-1} h_i, \qquad |\bar{\alpha}_k| \leq h_k \in \Lambda_k,$$

where $h_k$ can be any real number that satisfies the above inequality, which may or may not depend on the actual step-size $\bar{\alpha}_k$ computed by the boosting algorithm.

Using the definition of 1-norm for $\bar{f}$ and since $f_0 \in \text{span}(S)$, it is clear that for all $\varepsilon > 0$ we can choose a finite subset $S_k \subset S$, vector $\beta_k$ and vector $\bar{w}$ such that

$$\|\beta_k\|_1 = \sum_{g \in S_k} |\beta_k^g| \leq s_k + \varepsilon/2, \qquad \|\bar{w}\|_1 \leq \|\bar{f}\|_1 + \varepsilon/2.$$

It follows that with appropriate representation, the following inequality holds for all $\varepsilon > 0$:

$$(24) \qquad \Delta W_k \leq s_k + \|\bar{f}\|_1 + \varepsilon.$$

We now proceed to show that even in the worse case, the value $A(f_{k+1}) - A(\bar{f})$ decreases from $A(f_k) - A(\bar{f})$ by a reasonable quantity.

The basic idea is to upper bound the minimum of a set of numbers by an appropriately chosen weighted average of these numbers. This proof technique, which we shall call "averaging method," was used in [1, 20, 25, 37] for analysis of greedy-type algorithms.

For $h_k$ that satisfies (10), the symmetry of $\Lambda_k$ implies $h_k \, \text{sign}(\bar{w}^g - \beta_k^g) \in \Lambda_k$. Therefore the approximate minimization step (3) implies that for all $g \in S_k$, we have

$$A(f_{k+1}) \leq A(f_k + h_k s^g g) + \varepsilon_k, \qquad s^g = \text{sign}(\bar{w}^g - \beta_k^g).$$

Now multiply the above inequality by $|\bar{w}^g - \beta_k^g|$ and sum over $g \in S_k$; we obtain

$$(25) \qquad \Delta W_k(A(f_{k+1}) - \varepsilon_k) \leq \sum_{g \in S_k} |\beta_k^g - \bar{w}^g| A(f_k + h_k s^g g) =: B(h_k).$$

We only need to upper bound $B(h_k)$, which in turn gives an upper bound on $A(f_{k+1})$.

We recall a simple but important property of a convex function that follows directly from the definition of convexity of $A(f)$ as a function of $f$: for all $f_1, f_2$

$$(26) \qquad A(f_2) \geq A(f_1) + \nabla A(f_1)^T (f_2 - f_1).$$

If $A(f_k) - A(\bar{f}) < 0$, then $\Delta A(f_k) = 0$. From $0 \in \Lambda_k$ and (3), we obtain

$$A(f_{k+1}) - A(\bar{f}) \leq A(f_k) - A(\bar{f}) + \varepsilon_k \leq \bar{\varepsilon}_k,$$

which implies (13). Hence the lemma holds in this case. Therefore in the following, we assume that $A(f_k) - A(\bar{f}) \geq 0$.

Using Taylor expansion, we can bound each term on the right-hand side of (25) as

$$A(f_k + h_k s^g g) \leq A(f_k) + h_k s^g \nabla A(f_k)^T g + \frac{h_k^2}{2} \sup_{\xi \in [0,1]} A''_{f_k, g}(\xi h_k s^g).$$

Since Assumption 3.1 implies that

$$\sup_{\xi \in [0,1]} A''_{f_k, g}(\xi h_k s^g) = \sup_{\xi \in [0,1]} A''_{f_k + \xi h_k, g}(0) \leq M(\|f_k\|_1 + h_k),$$

we have

$$A(f_k + h_k s^g g) \leq A(f_k) + h_k s^g \nabla A(f_k)^T g + \frac{h_k^2}{2} M(\|f_k\|_1 + h_k).$$

Taking a weighted average, we have

$$B(h_k) = \sum_{g \in S_k} |\beta_k^g - \bar{w}^g| A(f_k + h_k s^g g)$$

$$\leq \sum_{g \in S_k} |\beta_k^g - \bar{w}^g| \left[ A(f_k) + \nabla A(f_k)^T h_k s^g g + \frac{h_k^2}{2} M(\|f_k\|_1 + h_k) \right]$$

$$= \Delta W_k A(f_k) + h_k \nabla A(f_k)^T (\bar{f} - f_k) + \frac{h_k^2}{2} \Delta W_k M(\|f_k\|_1 + h_k)$$

$$\leq \Delta W_k A(f_k) + h_k [A(\bar{f}) - A(f_k)] + \frac{h_k^2}{2} \Delta W_k M(\|f_k\|_1 + h_k).$$

The last inequality follows from (26). Now using (25) and the bound $\|f_k\|_1 + h_k \leq s_{k+1}$, we obtain

$$(A(f_{k+1}) - A(\bar{f})) - \varepsilon_k \leq \left( 1 - \frac{h_k}{\Delta W_k} \right)(A(f_k) - A(\bar{f})) + \frac{h_k^2}{2} M(s_{k+1}).$$

Now replace $\Delta W_k$ by the right-hand side of (24) with $\varepsilon \to 0$; we obtain the lemma.  $\square$

PROOF OF THE MULTISTEP ANALYSIS OR LEMMA 4.2. Note that for all $a \geq 0$,

$$\prod_{\ell=j}^{k} \left( 1 - \frac{h_\ell}{s_\ell + a} \right) = \exp\left( \sum_{\ell=j}^{k} \ln\left( 1 - \frac{s_{\ell+1} - s_\ell}{s_\ell + a} \right) \right)$$

$$\leq \exp\left( \sum_{\ell=j}^{k} -\frac{s_{\ell+1} - s_\ell}{s_\ell + a} \right) \leq \exp\left( -\int_{s_j}^{s_{k+1}} \frac{1}{v + a} \, dv \right)$$

$$= \frac{s_j + a}{s_{k+1} + a}.$$

By recursively applying (13) and using the above inequality, we obtain

$$\Delta A(f_{k+1}) \leq \prod_{\ell=0}^{k} \left(1 - \frac{h_\ell}{s_\ell + \|\bar{f}\|_1}\right) \Delta A(f_0) + \sum_{j=0}^{k} \prod_{\ell=j+1}^{k} \left(1 - \frac{h_\ell}{s_\ell + \|\bar{f}\|_1}\right) \bar{\varepsilon}_j$$

$$\leq \frac{s_0 + \|\bar{f}\|_1}{s_{k+1} + \|\bar{f}\|_1} \Delta A(f_0) + \sum_{j=0}^{k} \frac{s_{j+1} + \|\bar{f}\|_1}{s_{k+1} + \|\bar{f}\|_1} \bar{\varepsilon}_j. \qquad \square$$

**A.3. Discussion of step-size.** We have been deriving our results in the case of restricted step-size in which the crucial small step-size condition is explicit. In this section we investigate the case of unrestricted step-size under exact minimization, for which we show that the small step-size condition is actually implicit if the boosting algorithm converges. The implication is that the consistency (and rate of convergence) results can be extended to such a case, although the analysis becomes more complicated.

Let $\Lambda_k = R$ for all $k$, so that the size of $\bar{\alpha}_k$ in the boosting algorithm is unrestricted. For simplicity, we will only consider the case that $\sup_a M(a)$ is upper bounded by a constant $M$.

Interestingly enough, although the size of $\bar{\alpha}_k$ is not restricted in the boosting algorithm itself, for certain formulations the inequality $\sum_j \bar{\alpha}_j^2 < \infty$ still holds. Theorem 4.1 can then be applied to show the convergence of such boosting procedures. For convenience, we will impose the following additional assumption for the step-size $\bar{\alpha}_k$ in Algorithm 2.1:

$$(27) \qquad A(f_k + \bar{\alpha}_k \bar{g}_k) = \inf_{\alpha_k \in R} A(f_k + \alpha_k \bar{g}_k),$$

which means that given the selected basis function $\bar{g}_k$, the corresponding $\bar{\alpha}_k$ is chosen to be the exact minimizer.

LEMMA A.1. *Assume that $\bar{\alpha}_k$ satisfies* (27). *If there exists a positive constant $c$ such that*

$$\inf_k \inf_{\xi \in (0,1)} A''_{(1-\xi)f_k + \xi f_{k+1}, \bar{g}_k}(0) \geq c,$$

*then*

$$\sum_{j=0}^{k} \bar{\alpha}_j^2 \leq 2c^{-1}[A(f_0) - A(f_{k+1})].$$

PROOF. Since $\bar{\alpha}_k$ minimizes $A_{f_k, \bar{g}_k}(\alpha)$, $A'_{f_k, \bar{g}_k}(\bar{\alpha}_k) = 0$. Using Taylor expansion, we obtain

$$A_{f_k, \bar{g}_k}(0) = A_{f_k, \bar{g}_k}(\bar{\alpha}_k) + \tfrac{1}{2} A''_{f_k, \bar{g}_k}(\xi_k \bar{\alpha}_k) \bar{\alpha}_k^2,$$

where $\xi_k \in (0, 1)$. That is, $A(f_k) = A(f_{k+1}) + \frac{1}{2} A''_{f_k, \bar{g}_k}(\xi_k \bar{\alpha}_k) \bar{\alpha}_k^2$. By assumption, we have $A''_{f_k, \bar{g}_k}(\xi_k \bar{\alpha}_k) \geq c$. It follows that, $\forall j \geq 0$, $\bar{\alpha}_j^2 \leq 2c^{-1}[A(f_j) - A(f_{j+1})]$. We can obtain the lemma by summing from $j = 0$ to $k$.    $\square$

By combining Lemma A.1 and Corollary 4.1, we obtain:

COROLLARY A.1.    *Assume that $\sup_a M(a) < +\infty$ and $\varepsilon_j$ in (3) satisfies $\sum_{j=0}^{\infty} \varepsilon_j < \infty$. Assume also that in Algorithm 2.1 we let $\Lambda_k = R$ and let $\bar{\alpha}_k$ satisfy (27). If*

$$\inf_k \inf_{\xi \in (0,1)} A''_{(1-\xi)f_k + \xi f_{k+1}, \bar{g}_k}(0) > 0,$$

*then*

$$\lim_{k \to \infty} A(f_k) = \inf_{f \in \mathrm{span}(S)} A(f).$$

PROOF.    If $\lim_{k \to \infty} A(f_k) = -\infty$, then the conclusion is automatically true. Otherwise, Lemma A.1 implies that $\sum_{j=0}^{\infty} \bar{\alpha}_j^2 < \infty$. Now choose $h_j = |\bar{\alpha}_j| + 1/(j+1)$ in (10); we have $\sum_{j=0}^{\infty} h_j = \infty$, and $\sum_{j=0}^{\infty} h_j^2 < \infty$. The convergence now follows from Corollary 4.1.    $\square$

*Least squares loss.*    The convergence of unrestricted step-size boosting using least squares loss (matching pursuit) was studied in [29]. Since a scaling of the basis function does not change the algorithm, without loss of generality we can assume that $E_X g(X)^2 = 1$ for all $g \in S$ (assume $S$ does not contain function 0). In this case it is easy to check that for all $g \in S$,

$$A''_{f,g}(0) = E_X g(X)^2 = 1.$$

Therefore the conditions in Corollary A.1 are satisfied as long as $\sum_{j=0}^{\infty} \varepsilon_j < \infty$. This shows that the matching pursuit procedure converges, that is,

$$\lim_{k \to \infty} A(f_k) = \inf_{f \in \mathrm{span}(S)} A(f).$$

We would like to point out that for matching pursuit, the inequality in Lemma A.1 can be replaced by the equality

$$\sum_{j=0}^{k} \bar{\alpha}_j^2 = 2[A(f_0) - A(f_{k+1})],$$

which was referred to as "energy conservation" in [29], and was used there to prove the convergence.

*Exponential loss.* The convergence behavior of boosting with exponential loss was previously studied by Breiman [9] for $\pm 1$-trees under the assumption $\inf_x P(Y = 1|x)P(Y = -1|x) > 0$. Using exact computation, Breiman obtained an equality similar to the matching pursuit energy conservation equation. As part of the convergence analysis, the equality was used to show $\sum_{j=0}^{\infty} \bar{\alpha}_j^2 < \infty$.

The following lemma shows that under a more general condition, the convergence of unrestricted boosting with exponential loss follows directly from Corollary A.1. This result extends that of [9], but the condition still constrains the class of measures that generate the joint distribution of $(X, Y)$.

LEMMA A.2. *Assume that*

$$\inf_{g \in S} E_X |g(X)| \sqrt{P(Y = 1|X)P(Y = -1|X)} > 0.$$

*If $\bar{\alpha}_k$ satisfies* (27), *then $\inf_k \inf_{\xi \in (0,1)} A''_{(1-\xi)f_k + \xi f_{k+1}, \bar{g}_k}(0) > 0$. Hence $\sum_j \bar{\alpha}_j^2 < \infty$.*

PROOF. For notational simplicity, we let $q_{X,Y} = \exp(-f(X)Y)$. Recall that the direct computation of $A''_{f,g}(0)$ in Section A.1.2 yields

$$[E_{X,Y} q_{X,Y}]^2 A''_{f,g}(0)$$

$$= [E_{X,Y} g(X)^2 q_{X,Y}][E_{X,Y} q_{X,Y}] - [E_{X,Y} g(X)Y q_{X,Y}]^2$$

$$= [E_X g(X)^2 E_{Y|X} q_{X,Y}][E_X E_{Y|X} q_{X,Y}] - [E_X g(X) E_{Y|X} Y q_{X,Y}]^2$$

$$\geq [E_X g(X)^2 E_{Y|X} q_{X,Y}][E_X E_{Y|X} q_{X,Y}]$$

$$\quad - [E_X g(X)^2 |E_{Y|X} Y q_{X,Y}|][E_X |E_{Y|X} Y q_{X,Y}|]$$

$$\geq [E_X g(X)^2 E_Y q_{X,Y}] E_X [E_{Y|X} q_{X,Y} - |E_{Y|X} Y q_{X,Y}|]$$

$$\geq [E_X |g(X)| \sqrt{E_Y q_{X,Y}(E_{Y|X} q_{X,Y} - |E_{Y|X} Y q_{X,Y}|)}]^2$$

$$\geq [E_X |g(X)| \sqrt{2P(Y = 1|X)P(Y = -1|X)}]^2.$$

The first and the third inequalities follow from Cauchy–Schwarz, and the last inequality used the fact that $(a + b)((a + b) - |a - b|) \geq 2ab$. Now observe that $E_{X,Y} q_{X,Y} = \exp(A(f))$. The exact minimization (27) implies that $A(f_k) \leq A(f_0)$ for all $k \geq 0$. Therefore, using Jensen's inequality we know that $\forall \xi \in (0, 1), A((1 - \xi)f_k + \xi f_{k+1}) \leq A(f_0)$. This implies the desired inequality,

$$A''_{(1-\xi)f_k + \xi f_{k+1}, \bar{g}_k}(0)$$

$$\geq \exp(-2A(f_0))[E_X |\bar{g}_k(X)| \sqrt{2P(Y = 1|X)P(Y = -1|X)}]^2.$$

□

Although unrestricted step-size boosting procedures can be successful in certain cases, for general problems we are unable to prove convergence. In such cases the crucial condition of $\sum_{j=0}^{\infty} \bar{\alpha}_j^2 < \infty$, as required in the proof of Corollary A.1, can be violated. Although we do not have concrete examples at this point, we believe boosting may fail to converge when this condition is violated.

For example, for logistic regression we are unable to prove a result similar to Lemma A.2. The difficulty is caused by the near-linear behavior of the loss function toward negative infinity. This means that the second derivative is so small that we may take an extremely large step-size when $\bar{\alpha}_j$ is exactly minimized.

Intuitively, the difficulty associated with large $\bar{\alpha}_j$ is due to the potential problem of large oscillation in that a greedy step may search for a suboptimal direction, which needs to be corrected later on. If a large step is taken toward the suboptimal direction, then many more additional steps have to be taken to correct the mistake. If the additional steps are also large, then we may overcorrect and go to some other suboptimal directions. In general it becomes difficult to keep track of the overall effect.

**A.4. The relationship of AdaBoost and $L_1$-margin maximization.** Given a real-valued classification function $p(x)$, we consider the following discrete prediction rule:

$$
(28) \qquad y = \begin{cases} 1, & \text{if } p(x) \geq 0, \\ -1, & \text{if } p(x) < 0. \end{cases}
$$

Its classification error [for simplicity we ignore the point $p(x) = 0$, which is assumed to occur rarely] is given by

$$
L_\gamma(p(x), y) = \begin{cases} 1, & \text{if } p(x)y \leq \gamma, \\ 0, & \text{if } p(x)y > \gamma, \end{cases}
$$

with $\gamma = 0$. In general, we may consider $\gamma \geq 0$ and the parameter $\gamma \geq 0$ is often referred to as margin, and we shall call the corresponding error function $L_\gamma$ margin error.

In [33] the authors proved that under appropriate assumptions on the base learner, the expected margin error $L_\gamma$ with a positive margin $\gamma > 0$ also decreases exponentially. It follows that regularity assumptions of weak learning for AdaBoost imply the following margin condition: there exists $\gamma > 0$ such that $\inf_{f \in \text{span}(S), \|f\|_1 = 1} L_\gamma(f, y) = 0$, which in turn implies the inequality for all $s > 0$,

$$
(29) \qquad \inf_{f \in \text{span}(S), \|f\|_1 = 1} E_{X,Y} \exp(-sf(X)Y) \leq \exp(-\gamma s).
$$

We now show that under (29) the expected margin errors (with small margin) from Algorithm 2.1 may decrease exponentially. A similar analysis was given in [37]. However, the boosting procedure considered there was modified so that the estimator always stays in the scaled convex hull of the basis functions. This restriction is removed in the current analysis:

$$f_0 = 0, \qquad \sup \Lambda_k \le h_k, \qquad \varepsilon_k \le h_k^2/2.$$

Note that this implies that $\bar{\varepsilon}_k \le h_k^2$ for all $k$.

Now applying (15) with $\tilde{f} = sf$ for any $s > 0$ and letting $f$ approach the minimum in (29), we obtain (recall $\|f\|_1 = 1$)

$$A(f_k) \le -s\gamma \frac{s_k}{s_k + s} + \sum_{j=1}^{k} \frac{s_j + s}{s_k + s} \bar{\varepsilon}_{j-1} \le -s\gamma \frac{s_k}{s_k + s} + \sum_{j=0}^{k-1} h_j^2.$$

Now let $s \to \infty$; we have

$$A(f_k) \le -\gamma s_k + \sum_{j=0}^{k-1} h_j^2.$$

Assume we pick a constant $h < \gamma$ and let $h_k = h$; then

$$(30) \qquad E_{X,Y} \exp(-f_k(X)Y) \le \exp(-kh(\gamma - h)),$$

which implies that the margin error decreases exponentially for all margins less than $\gamma - h$. To see this, consider $\gamma' < \gamma - h$. Since $\|f_k\|_1 \le kh$, we have from (30),

$$L_{\gamma'}(f_k(x)/\|f_k\|_1, y) \le P(f_k(X)Y \le kh\gamma')$$
$$\le E_{X,Y} \exp(-f_k(X)Y + kh\gamma') \le \exp(-kh(\gamma - h - \gamma')).$$

Therefore

$$\lim_{k \to \infty} L_{\gamma'}(f_k(x)/\|f_k\|_1, y) = 0.$$

This implies that as $h \to 0$, $f_k(x)/\|f_k\|_1$ achieves a margin that is within $h$ of the maximum possible. Therefore, when $h \to$ and $k \to \infty$, $f_k(x)/\|f_k\|_1$ approaches a maximum margin separator.

Note that in this particular case we allow a small step-size ($h < \gamma$), which violates the condition $\sum_k h_k^2 < \infty$ imposed for the boosting algorithm to converge. However, this condition that prevents large oscillation from occurring is only a sufficient condition to guarantee convergence. For specific problems, especially when $\inf_{f \in \text{span}(S)} A(f) = -\infty$, it is still possible to achieve convergence even if the condition is violated.

## REFERENCES

[1] BARRON, A. (1993). Universal approximation bounds for superpositions of a sigmoidal function. *IEEE Trans. Inform. Theory* **39** 930–945. MR1237720

[2] BARTLETT, P. L., BOUSQUET, O. and MENDELSON, S. (2005). Local Rademacher complexities. *Ann. Statist.* **33** 1497–1537.

[3] BARTLETT, P. L., JORDAN, M. and MCAULIFFE, J. (2005). Convexity, classification, and risk bounds. *J. Amer. Statist. Assoc.* To appear.

[4] BARTLETT, P. L. and MENDELSON, S. (2002). Rademacher and Gaussian complexities: Risk bounds and structural results. *J. Mach. Learn. Res.* **3** 463–482. MR1984026

[5] BLANCHARD, G., LUGOSI, G. and VAYATIS, N. (2004). On the rate of convergence of regularized boosting classifiers. *J. Mach. Learn. Res.* **4** 861–894. MR2076000

[6] BOUSQUET, O., KOLTCHINSKII, V. and PANCHENKO, D. (2002). Some local measures of complexity of convex hulls and generalization bounds. *Computational Learning Theory. Lecture Notes in Artificial Intelligence* **2375** 59–73. Springer, Berlin. MR2040405

[7] BREIMAN, L. (1998). Arcing classifiers (with discussion). *Ann. Statist.* **26** 801–849. MR1635406

[8] BREIMAN, L. (1999). Prediction games and arcing algorithms. *Neural Computation* **11** 1493–1517.

[9] BREIMAN, L. (2004). Population theory for boosting ensembles. *Ann. Statist.* **32** 1–11. MR2050998

[10] BÜHLMANN, P. (2002). Consistency for $L_2$ boosting and matching pursuit with trees and tree-type basis functions. Technical report, ETH Zürich.

[11] BÜHLMANN, P. and YU, B. (2003). Boosting with the $L_2$ loss: Regression and classification. *J. Amer. Statist. Assoc.* **98** 324–339. MR1995709

[12] COLLINS, M., SCHAPIRE, R. E. and SINGER, Y. (2002). Logistic regression, AdaBoost and Bregman distances. *Machine Learning* **48** 253–285.

[13] FREUND, Y. and SCHAPIRE, R. (1997). A decision-theoretic generalization of on-line learning and an application to boosting. *J. Comput. System Sci.* **55** 119–139. MR1473055

[14] FRIEDMAN, J. H. (2001). Greedy function approximation: A gradient boosting machine. *Ann. Statist.* **29** 1189–1232. MR1873328

[15] FRIEDMAN, J. H., HASTIE, T. and TIBSHIRANI, R. (2000). Additive logistic regression: A statistical view of boosting (with discussion). *Ann. Statist.* **28** 337–407. MR1790002

[16] GROVE, A. and SCHUURMANS, D. (1998). Boosting in the limit: Maximizing the margin of learned ensembles. In *Proc. Fifteenth National Conference on Artificial Intelligence* 692–699. AAAI Press, Menlo Park, CA.

[17] HASTIE, T. J. and TIBSHIRANI, R. J. (1990). *Generalized Additive Models.* Chapman and Hall, London. MR1082147

[18] HASTIE, T., TIBSHIRANI, R. and FRIEDMAN, J. H. (2001). *The Elements of Statistical Learning.* Springer, New York. MR1851606

[19] JIANG, W. (2004). Process consistency for AdaBoost. *Ann. Statist.* **32** 13–29. MR2050999

[20] JONES, L. (1992). A simple lemma on greedy approximation in Hilbert space and convergence rates for projection pursuit regression and neural network training. *Ann. Statist.* **20** 608–613. MR1150368

[21] KOLTCHINSKII, V. and PANCHENKO, D. (2002). Empirical margin distributions and bounding the generalization error of combined classifiers. *Ann. Statist.* **30** 1–50. MR1892654

[22] KOLTCHINSKII, V. and PANCHENKO, D. (2005). Complexities of convex combinations and bounding the generalization error in classification. *Ann. Statist.* **33** 1455–1496.

[23] KOLTCHINSKII, V., PANCHENKO, D. and LOZANO, F. (2001). Further explanation of the effectiveness of voting methods: The game between margins and weights. *Computational Learning Theory. Lecture Notes in Artificial Intelligence* **2111** 241–255. Springer, Berlin. MR2042039

[24] LEDOUX, M. and TALAGRAND, M. (1991). *Probability in Banach Spaces: Isoperimetry and Processes.* Springer, Berlin. MR1122615

[25] LEE, W., BARTLETT, P. and WILLIAMSON, R. (1996). Efficient agnostic learning of neural networks with bounded fan-in. *IEEE Trans. Inform. Theory* **42** 2118–2132. MR1447518

[26] LESHNO, M., LIN, Y. V., PINKUS, A. and SCHOCKEN, S. (1993). Multilayer feedforward networks with a non-polynomial activation function can approximate any function. *Neural Networks* **6** 861–867.

[27] LI, F. and YANG, Y. (2003). A loss function analysis for classification methods in text categorization. In *Proc. 20th International Conference on Machine Learning* **2** 472–479. AAAI Press, Menlo Park, CA.

[28] LUGOSI, G. and VAYATIS, N. (2004). On the Bayes-risk consistency of regularized boosting methods. *Ann. Statist.* **32** 30–55. MR2051000

[29] MALLAT, S. and ZHANG, Z. (1993). Matching pursuits with time-frequency dictionaries. *IEEE Trans. Signal Process.* **41** 3397–3415.

[30] MANNOR, S., MEIR, R. and ZHANG, T. (2003). Greedy algorithms for classification—consistency, convergence rates, and adaptivity. *J. Mach. Learn. Res.* **4** 713–742. MR2072266

[31] MASON, L., BAXTER, J., BARTLETT, P. L. and FREAN, M. (2000). Functional gradient techniques for combining hypotheses. In *Advances in Large Margin Classifiers* (A. J. Smola, P. L. Bartlett, B. Schölkopf and D. Schuurmans, eds.) 221–246. MIT Press. MR1820960

[32] MEIR, R. and ZHANG, T. (2003). Generalization error bounds for Bayesian mixture algorithms. *J. Mach. Learn. Res.* **4** 839–860.

[33] SCHAPIRE, R., FREUND, Y., BARTLETT, P. L. and LEE, W. (1998). Boosting the margin: A new explanation for the effectiveness of voting methods. *Ann. Statist.* **26** 1651–1686. MR1673273

[34] SCHAPIRE, R. and SINGER, Y. (1999). Improved boosting algorithms using confidence-rated predictions. *Machine Learning* **37** 297–336. MR1811573

[35] VAN DER VAART, A. W. and WELLNER, J. A. (1996). *Weak Convergence and Empirical Processes. With Applications to Statistics.* Springer, New York. MR1385671

[36] VAPNIK, V. (1998). *Statistical Learning Theory.* Wiley, New York. MR1641250

[37] ZHANG, T. (2003). Sequential greedy approximation for certain convex optimization problems. *IEEE Trans. Inform. Theory* **49** 682–691. MR1967192

[38] ZHANG, T. (2004). Statistical behavior and consistency of classification methods based on convex risk minimization. *Ann. Statist.* **32** 56–85. MR2051051

[39] ZHANG, T. and OLES, F. J. (2001). Text categorization based on regularized linear classification methods. *Information Retrieval* **4** 5–31.

IBM T. J. WATSON RESEARCH CENTER
YORKTOWN HEIGHTS, NEW YORK 10598
USA
E-MAIL: tzhang@watson.ibm.com

DEPARTMENT OF STATISTICS
UNIVERSITY OF CALIFORNIA
BERKELEY, CALIFORNIA 94720-3860
USA
E-MAIL: binyu@stat.berkeley.edu